\documentclass[11pt]{article}
\usepackage{amsfonts}
\usepackage{latexsym,amsfonts,amssymb}
\usepackage{indentfirst}
\usepackage{amsmath}
\usepackage{mathrsfs}
\usepackage{color}

\title{Sharp well-posedness and ill-posedness in Fourier-Besov spaces for the viscous primitive equations of geophysics\thanks{Project supported
by National Natural Science Foundation of China under the grant
number 11571381.}}

\author{Jinyi Sun\thanks{Corresponding author.  E-mails: sunjinyi333@163.com (J. Sun) and cuishb@mail.sysu.edu.cn (S. Cui).}
\ \ and \ \ Shangbin Cui\\
{\small Department of Mathematics, Sun Yat-Sen University, Guangzhou 510275, P. R. China}}

\date{}

\begin{document}
\maketitle

\begin{abstract}
We study well-posedness and ill-posedness for Cauchy problem of the three-dimensional viscous primitive equations
describing the large scale ocean and atmosphere dynamics. By using the Littlewood-Paley analysis technique, in particular
Chemin-Lerner's localization method, we prove that the Cauchy problem with Prandtl number $P=1$ is locally
well-posed in the Fourier-Besov spaces $[\dot{FB}^{2-\frac{3}{p}}_{p,r}(\mathbb{R}^3)]^4$ for $1<p\leq\infty,1\leq r<\infty$ and
$[\dot{FB}^{-1}_{1,r}(\mathbb{R}^3)]^4$ for $1\leq r\leq 2$, and globally well-posed in these spaces when the initial data $(u_0,\theta_0)$
are small. We also prove that such problem is ill-posed in $[\dot{FB}^{-1}_{1,r}(\mathbb{R}^3)]^4$ for $2<r\leq\infty$,
showing that the results stated above are sharp.



\end{abstract}

\section{Introduction}

The viscous primitive equations are a fundamental mathematical model in the field of fluid geophysics. It describes the large scale ocean and atmosphere dynamics, cf. the monograph \cite{Cushman, Pedlosky, Salmon}, for instance.
The model reads as follows:
$$
\left\{
  \begin{array}{ll}
   \partial_t u-\nu\Delta u+\Omega e_3\times u+(u\cdot\nabla)u+\nabla p=g\theta e_3~~~~~& \textrm{in}~~\mathbb{R}^3\times (0,\infty), \\
   \partial_t \theta-\mu\Delta \theta+(u\cdot\nabla)\theta=-\mathcal{N}^2u_3 & \textrm{in}~~\mathbb{R}^3\times (0,\infty), \\
   \textrm{div} u=0 & \textrm{in}~~\mathbb{R}^3\times (0,\infty), 
  \end{array}
\right.
\eqno(1.1)
$$
where the unknown functions $u=(u_1, u_2, u_3)$ and $p$ denote the fluid velocity and the pressure, respectively, and $\theta$ is a scalar function representing the density fluctuation in the fluid (in the case of the ocean it depends on the temperature and the salinity, and in the case of the atmosphere it depends on the temperature), and $\nu$, $\mu$ and $g$ are positive constants related to viscosity, diffusivity and gravity, respectively. Moreover, $\Omega$ is the so-called Coriolis parameter, a real constant which is twice the angular velocity of the rotation around the vertical unit vector $e_3=(0,0,1)$, and $\mathcal{N}$ is the stratification parameter, a nonnegative constant representing the Brunt-V\"{a}is\"{a}l\"{a} wave frequency.
The ratio $P:=\frac{\nu}{\mu}$ is known as the Prandtl number and $B:=\frac{\Omega}{\mathcal{N}}$ is essentially the ``Burger'' number of geophysics. We refer the reader to see
\cite{Babin 1999, Cushman, Pedlosky, Salmon} for derivation of this model and more detailed discussions on its physical background.

If $\theta\equiv0$, $\mathcal{N}=0$ and $\Omega=0$, then (1.1) reduces to the classical incompressible Navier-Stokes equations
$$
\left\{
  \begin{array}{ll}
   \partial_t u-\nu\Delta u+(u\cdot\nabla)u=-\nabla p~~~~~ & \textrm{in}~~\mathbb{R}^3\times (0,\infty), \\
   \textrm{div} u=0 & \textrm{in}~~\mathbb{R}^3\times (0,\infty), 
  \end{array}
\right.
\eqno(1.2)
$$
which has drawn great attention during the past fifty more years. 
In 1964 Fujita and Kato \cite{Fujita} obtained the first result on well-posedness of the Cauchy problem of (1.2) and proved that it is locally well-posed in $H^s(\mathbb{R}^3)$ for $s\geq \frac{1}{2}$ and globally well-posed in $H^\frac{1}{2}(\mathbb{R}^3)$ for small initial data. These results were later extended to various other function spaces, cf. \cite{Cannone 1997, Cannone 2012, Cui,Giga 1986,Giga 1989,Kato,Koch,Lei,Planchon} and references therein. Particularly worth mentioning is that well-posedness has been established in $\dot{B}^{-1+\frac{3}{p}}_{p, r}(\mathbb{R}^3)$ for $1\leq p<\infty$, $1\leq r\leq \infty$ by Cannone \cite{Cannone 1997} and Planchon \cite{Planchon}, in $BMO^{-1}$ by Koch and Tataru \cite{Koch} and in $B^{-1,\sigma}_{\infty,q}(\mathbb{R}^3)$ for $1\leq q\leq\infty$ and $\sigma\geq 1-\min\{1-\frac{1}{q},\frac{1}{q}\}$ by Cui \cite{Cui}.
On the other hand, the ill-posedness for Cauchy problem of (1.2) has been proved in $\dot{B}^{-1}_{\infty, q}(\mathbb{R}^3)$ for $1\leq q\leq \infty$,
in $\dot F^{-1}_{\infty, r}(\mathbb{R}^3)$ for $2<r\leq\infty$ and in $B^{-1,\sigma}_{\infty,q}(\mathbb{R}^3)$ for $1\leq q\leq\infty$ and $0\leq \sigma<1-\min\{1-\frac{1}{q},\frac{1}{q}\}$, cf. Bourgain and Pavlovi\'{c} \cite{Bourgain}, Yoneda
\cite{Yoneda}, Wang \cite{Wang} and Cui \cite{Cui}, which imply that the well-posedness results obtained in \cite{Cannone 1997,Planchon,Koch,Cui} are sharp.

If only $\theta\equiv0$ and $\mathcal{N}=0$ but $\Omega\neq 0$, then (1.1) reduces to the
incompressible rotating Navier-Stokes equation
$$
\left\{
  \begin{array}{ll}
   \partial_t u-\nu\Delta u+\Omega e_3\times u+(u\cdot\nabla)u+\nabla p=0\ \ \ \ \ & \textrm{in}~~\mathbb{R}^3\times (0,\infty), \\
   \textrm{div} u=0 & \textrm{in}~~\mathbb{R}^3\times (0,\infty). 
  \end{array}
\right.
\eqno(1.3)
$$
The topic of well-posedness and ill-posedness of the Cauchy problem of (1.3) have also been widely studied, cf. \cite{Babin 1999,
Chemin 2002, Chemin 2006, Fang, Fang 2015, Giga, Hieber, Iwabuchi 2013, Iwabuchi 2014, Koh, Konieczny, Sun} and the references therein. In particular, it has been
proved that the Cauchy problem of (1.3) is globally well-posed for small initial data in $FM^{-1}_0(\mathbb{R}^3)$ by Giga et al. \cite{Giga}, in $H^{\frac{1}{2}}(\mathbb{R}^3)$ by Hieber et al. \cite{Hieber}, in $\dot{FB}^{2-\frac{3}{p}}_{p,p}(\mathbb{R}^3)$ for $3<p<\infty$ and $\dot{FB}^{-1}_{1,1}(\mathbb{R}^3)\cap
\dot{FB}^{0}_{1,1}(\mathbb{R}^3)$ by Konieczny and Yoneda \cite{Konieczny}, in $\dot{FB}^{2-\frac{3}{p}}_{p,r}(\mathbb{R}^3)$ for $1<p\leq\infty$ and $1\leq r<\infty$ by Fang et al. \cite{Fang} and in $\dot{FB}^{-1}_{1,2}(\mathbb{R}^3)$ by Iwabuchi and Takada
\cite{Iwabuchi 2014}. Moreover, in \cite{Iwabuchi 2014}, the ill-posedness has been verified in $\dot{FB}^{-1}_{1,r}(\mathbb{R}^3)$ for $2<r\leq\infty$, which implies that  the well-posedness results in $\dot{FB}^{2-\frac{3}{p}}_{p,r}(\mathbb{R}^3)$
for $1<p\leq\infty$, $1\leq r<\infty$ and $\dot{FB}^{-1}_{1,2}(\mathbb{R}^3)$ obtained in \cite{Fang} and \cite{Iwabuchi 2014} are sharp.

In this paper we study well-posedness and ill-posedness of the Cauchy problem of (1.1):
$$
\left\{
  \begin{array}{ll}
   \partial_t u-\nu\Delta u+\Omega e_3\times u+(u\cdot\nabla)u+\nabla p=g\theta e_3~~~&\textrm{in}~~\mathbb{R}^3\times (0,\infty), \\
   \partial_t \theta-\mu\Delta \theta+(u\cdot\nabla)\theta=-\mathcal{N}^2u_3 & \textrm{in}~~\mathbb{R}^3\times (0,\infty), \\
   \textrm{div} u=0  & \textrm{in}~~\mathbb{R}^3\times (0,\infty), \\
   u|_{t=0}=u_0,~~\theta|_{t=0}=\theta_0 & \textrm{in}~~\mathbb{R}^3,
  \end{array}
\right.
\eqno(1.4)
$$
where $u_0$ and $\theta_0$ are given initial functions. A short review of existing work on this topic is as follows.
In \cite{Babin 1999}, by taking full advantage of the absence of
resonances between the fast rotation and the nonlinear advection, Babin, Maholov and Nicolaenko obtained global well-posedness of problem (1.4) in $[H^s(\mathbb{T}^3)]^4$ with $s\geq3/4$ for small initial data when the stratification parameter $\mathcal{N}$
is sufficiently large. Later on, by constructing the solution of a quasi-geostrophic system related to
equations (1.1) and using some Strichartz-type estimates, Charve \cite{Charve 2004} verified global well-posedness of problem (1.4) in $[\dot H^{\frac{1}{2}}(\mathbb{R}^3)\cap
\dot H^{1}(\mathbb{R}^3)]^4$ for arbitrary (i.e., not necessarily small) initial data under the assumptions that both $\Omega$ and $\mathcal{N}$ are
sufficiently large (depending on the scale of the initial data). In \cite{Charve 2008}, Charve further considered global well-posedness of (1.4) in
less regular initial value spaces. We also mention the interesting work of Charve and Ngo \cite{Charve 2011} on well-posedness of the problem
(1.4) with anisotropic viscosities. Recently, Koba, Mahalov and Yoneda \cite{Koba}
proved global well-posedness of problem (1.4) for any given $(u_0,\theta_0)
\in[\dot H^{\frac{1}{2}}(\mathbb{R}^3)\cap\dot H^{1}(\mathbb{R}^3)]^4$ with $\partial_2u_0^1-\partial_1u_0^2=0$ in the special case Prandtl number $P=1$ provided either condition $(a)$: the absolute value of ``Burger'' number $|B|<\sqrt{g}$ and $\mathcal{N}$ is sufficiently large
(depending on the scale of the initial data) or condition $(b)$: the absolute value of ``Burger'' number $|B|>\sqrt{g}$ and both $\Omega$ and $\mathcal{N}$ are sufficiently large (depending on the scale
of initial data) holds. They also obtained the global well-posedness of problem (1.4) for uniformly small initial data in $[\dot{H}^{\frac{1}{2}}(\mathbb{R}^3)]^4$.
For the other related studies on the viscous primitive equations (1.1), we refer the interested reader to \cite{Babin 1997,Babin 2000,Babin 2001,Cao 2007, Cao 2012, Cao 2015, Charve 2005,Charve 2006,Charve 2014,Chemin 1997}.

From the above review we see that concerning the well-posedness issue, there is a big gap between the standard Navier-Stokes initial value problem (1.2) and the initial value problem (1.4) of the primitive equations: We know that the problem (1.2) is well-posed in a lot of function spaces
of negative regularity indices, so that it has at least local-in-time solutions for a large group of rough initial data. For the problem (1.4), however, existing results only show that it is locally well-posed in some function spaces of the regularity index $s\geq\frac{1}{2}$. It is natural to ask whether similar results concerning well-posedness and ill-posedness of the initial value problem in function spaces of negative regularity indices as for the Navier-Stokes equations can be established for the primitive equations. This is the main motivation of the present work. However, due to the influence of the oscillations caused by the rotation (i.e., the term $\Omega e_3\times u$) and the stratification (i.e., the terms $g\theta e_3$ and $\mathcal{N}^2u_3$), a big portion of the integral estimates for the Stokes semigroup $\{\textrm{e}^{t\Delta}\}_{t\geq 0}$ (which relates to the Navier-Stokes equations) do not work for the Stokes-Coriolis-Stratification semigroup $\{T_{\Omega,N}(t)\}_{t\geq 0}$ (see Section 2.1 for the definition) related to the primitive equations. Consequently, the usual function spaces used in the study of the Navier-Stokes equations such as the homogeneous and inhomogeneous Besov spaces $\dot{B}^s_{pr}(\mathbb{R}^3)$ and $B^s_{pr}(\mathbb{R}^3)$ and the space $BMO^{-1}(\mathbb{R}^3)$ are not suitable for the primitive equations.
In this work, as in \cite{Fang,Iwabuchi 2014,Konieczny}, we use the Fourier-Besov spaces $[\dot{FB}^{2-\frac{3}{p}}_{p,r}(\mathbb{R}^3)]^4$ ($1\leq p,r\leq\infty$) as the initial value space to study well-posedness issue of the primitive equations.

To give the precise statements of our main results, we first recall the definitions of the homogeneous Besov spaces $\dot{B}^{s}_{p,r}(\mathbb{R}^3)$ and the Fourier-Besov space $\dot{FB}^s_{p,r}(\mathbb{R}^3)$.
As usual we denote by $\mathscr{S}(\mathbb{R}^3)$
the space of Schwartz functions on $\mathbb{R}^3$, and by $\mathscr{S}'(\mathbb{R}^3)$ the space of tempered distributions on $\mathbb{R}^3$.
Choose two radial function $\varphi, \psi\in\mathscr{S}(\mathbb{R}^3)$ such that their Fourier transform $\hat{\varphi}$ and $\hat{\psi}$ satisfies the following properties:
$$\textrm{supp}~\hat{\varphi}\subset \mathcal{B}:=\{\xi\in
\mathbb{R}^3:|\xi|\leq\frac{4}{3}\},$$
$$\textrm{supp}~\hat{\psi}\subset \mathcal{C}:=\{\xi\in \mathbb{R}^3:\frac{3}{4}\leq|\xi|\leq\frac{8}{3}\},
$$
and, furthermore,
$$\sum_{j\in\mathbb{Z}} \hat{\psi}(2^{-j}\xi)=1 ~\quad \textrm{for all } \xi\in \mathbb{R}^3\setminus\{0\}.$$
Let $\varphi_j(x):=2^{3j}\varphi(2^jx)$ and $\psi_j(x):=2^{3j}\psi(2^jx)$ for all $j\in\mathbb{Z}$. We define by $\Delta_j$ and $S_j$ the following operators in
$\mathscr{S}'(\mathbb{R}^3)$:
$$\Delta_j f:=\psi_j \ast f~~\textrm{and}~~S_j f:=\varphi_j\ast f  \quad \textrm{for}~j\in \mathbb{Z}\ \ \text{and}\ \ f\in \mathscr{S}'(\mathbb{R}^3).$$
Define $\mathscr{S}'_h(\mathbb{R}^3):=\mathscr{S}'(\mathbb{R}^3)/\mathcal{P}[\mathbb{R}^3]$, where $\mathcal{P}[\mathbb{R}^3]$ denotes the linear space of polynomials on $\mathbb{R}^3$. It is known that  there hold the following decompositions
$$f=\sum_{j\in\mathbb{Z}}\Delta_j f \ \ \ \mbox{\rm and}\ \ \ \
  S_j f=\sum_{j'\leq j-1}\Delta_{j'}f  \ \ \ \mbox{\rm in}\ \ \ \mathscr{S}'_h(\mathbb{R}^3),
$$
see \cite{Bahouri} for reference. With our choice of $\varphi$ and $\psi$, it is easy to verify that
$$\Delta_j\Delta_k f=0\ \ \ \ \ \hbox{if}\ |j-k|\geq 2 \ \hbox{and}$$
$$\Delta_j(S_{k-1}f\Delta_kf)=0\ \ \ \ \hbox{if}\ |j-k|\geq 5.$$


{\bf Definition 1.1}\ \ (\cite{Bahouri}) \ \ For $s\in\mathbb{R}$ and $1\leq p,r\leq \infty$, the homogeneous Besov spaces $\dot{B}^s_{p,r}(\mathbb{R}^3)$ consists of those distributions $u$ in
$\mathscr{S}'_h(\mathbb{R}^3)$ such that
$$\|u\|_{\dot{B}^{s}_{p,r}}:=\Big\|\Big\{2^{js}\|
\Delta_ju\|_{L^p}\Big\}_{j\in\mathbb{Z}}\Big\|_{\ell^r(\mathbb{Z})}<\infty.$$

{\bf Definition 1.2}\ \ (\cite{Iwabuchi 2014, Konieczny}) \ \ For $s\in\mathbb{R}$ and $1\leq p,r\leq \infty$, the Fourier-Besov space $\dot{FB}^s_{p,r}(\mathbb{R}^3)$ is defined to be the
set of all tempered distributions $u\in \mathscr{S}'_h(\mathbb{R}^3) $ such that
$$\|u\|_{\dot{FB}^{s}_{p,r}}:=\Big\|\Big\{2^{js}\|
\widehat{\Delta_ju}\|_{L^p}\Big\}_{j\in\mathbb{Z}}\Big\|_{\ell^r(\mathbb{Z})}<\infty.$$

Fourier-Besov spaces are introduced in the literature very recently, and sometime they are entitled with other different names. An early paper by
Cannone and Karch \cite{Cannone 2004} studied well-posedness of the Cauchy problem of (1.2) in the space $\mathscr{PM}^2$, which is in fact the space
$\dot {FB}^2_{\infty,\infty}(\mathbb{R}^3)$. Iwabuchi \cite{Iwabuchi 2011} investigated well-posedness and ill-posedness for Cauchy problem of
Keller-Segel system in $\mathcal{\dot B}^{-2}_{q} (2\leq q\leq \infty)$, which is in fact the space $\dot {FB}^{-2}_{1,q}(\mathbb{R}^3)(2\leq q\leq
\infty)$. Lei and Lin \cite{Lei} proved global existence of mild solutions to the Cauchy problem of (1.2) in $\mathcal{X}^{-1}$, which is in fact
equal to the space $\dot {FB}^{-1}_{1,1}(\mathbb{R}^3)$. Cannone and Wu \cite{Cannone 2012} extended the result in \cite{Lei} to the Fourier-Herz
space $\mathscr{\dot B}^{-1}_q$ $(1\leq q\leq 2)$, which is in fact the space $\dot {FB}^{-1}_{1,q}(\mathbb{R}^3)$ $(1\leq q\leq 2)$, and Liu and
Zhao \cite{Liu} considered the global well-posedness of the Cauchy problem of generalized magneto-hydrodynamic equations in the Fourier-Herz spaces
$\mathscr{\dot B}^{-(2\beta-1)}_q$ $(1\leq q\leq 2)$, which is in fact the space $\dot {FB}^{-(2\beta-1)}_{1,q}(\mathbb{R}^3)$ $(1\leq q\leq 2)$.
Systematic utilization of the Fourier-Besov spaces $\dot{FB}^{2-\frac{3}{p}}_{p,r}(\mathbb{R}^3)$ first appeared in the references \cite{Iwabuchi 2014, Konieczny} mentioned above.

Note that the definition of the Fourier-Besov spaces $\dot{FB}^s_{p,r}(\mathbb{R}^3)$ has some similar feature with that of the classical homogeneous
Besov spaces $\dot{B}^s_{p,r}(\mathbb{R}^3)$: they both measure regularity of a function $u$ with index $s$ which depicts the decay or increment
speed of its Fourier transform $\hat{u}(\xi)$ as $\xi\to\infty$ via dyadic decomposition. The only difference in their definitions is that for the
Fourier-Besov spaces such measurement is made purely in the frequency space, while for the Besov spaces $\dot{B}^s_{p,r}(\mathbb{R}^3)$ this is done in both frequency and physical
spaces jointly. Precisely because in the definition of the Fourier-Besov spaces regularity of a function is measured purely through its frequencies, they are
very useful in the study of partial differential equations which are not of purely dissipative type but instead of dissipative and dispersive joint
type, such as the equations (1.1) and (1.3). Relations between the Fourier-Besov spaces $\dot{FB}^s_{p,r}(\mathbb{R}^3)$ and the homogeneous Besov
spaces $\dot{B}^s_{p,r}(\mathbb{R}^3)$ are as follows:
$$
  \dot{FB}^s_{2,r}(\mathbb{R}^3)=\dot{B}^s_{2,r}(\mathbb{R}^3)\ \ \ \text{for}\;\; s\in\mathbb{R}\;\, \text{and}\;\, r\in[1,\infty]
$$
and
$$
  \dot{FB}^s_{p,r}(\mathbb{R}^3)\hookrightarrow\dot{B}^s_{p',r}(\mathbb{R}^3)\ \ \ \ \text{and}\ \ \
  \dot{B}^s_{p,r}(\mathbb{R}^3)\hookrightarrow\dot{FB}^s_{p',r}(\mathbb{R}^3)
$$
for $s\in\mathbb{R},~p\in[1,2]$ and $r\in[1,\infty]$. These relations can be easily proved by using the Plancherel identity and the Hausdorff-Young
inequality, cf. \cite{Iwabuchi 2014, Konieczny}.

\medskip

{\bf Definition 1.3}\ \  For $T>0$, $s\in\mathbb{R}$ and $1\leq r,\delta\leq \infty$, the Chemin-Lerner type space $\tilde{L}^\delta(0,T;
\dot{FB}^{s}_{p,r}(\mathbb{R}^3))$ built on $\dot{FB}^s_{p,r}(\mathbb{R}^3)$ is defined to be the set of all strongly measurable functions $u:(0,T)\to
\dot{FB}^s_{p,r}(\mathbb{R}^3)$ such that
$$\|u\|_{\tilde{L}^\delta(0,T; \dot{FB}^{s}_{p,r})}:=\Big\|\Big\{2^{js}\|
\widehat{\Delta_ju}\|_{L^\delta(0,T; L^p)}\Big\}_{j\in\mathbb{Z}}\Big\|_{l^r(\mathbb{Z})}<\infty.$$

\medskip

The main results of this paper are the following three theorems:

\medskip

{\bf Theorem 1.4}\, \, Let Prandtl number $P=1$, i.e., $\nu=\mu$. Assume that $p\in(1,\infty]$ and $r\in[1, \infty)$. Then for any $(u_0, \theta_0)\in
\big[\dot{FB}^{2-\frac{3}{p}}_{p,r}(\mathbb{R}^3)\big]^{4}$ satisfying $\textrm{div} u_0=0$, there exists corresponding $T>0$
such that problem (1.4) possesses a unique mild solution $(u, \theta)\in \big[C\big([0,T],\dot{FB}^{2-\frac{3}{p}}_{p,r}(\mathbb{R}^3)\big)\big]^{4}\cap
X_T^\alpha$, where
$$X_T^\alpha:=[\tilde{L}^{\frac{2}{1+\alpha}}(0,T;\dot{FB}^{3-\frac{3}{p}+\alpha}_{p,r}(\mathbb{R}^3))]^4
\cap [\tilde{L}^{\frac{2}{1-\alpha}}(0,T;\dot{FB}^{3-\frac{3}{p}-\alpha}_{p,r}(\mathbb{R}^3))]^4$$
(for an arbitrarily chosen but fixed number $\alpha\in(0,1)$).
Moreover, there exists a constant $C>0$ independent of $\nu$,
$\Omega$ and $\mathcal{N}$ such that if
$$\|(u_0, \sqrt{g}\theta_0/\mathcal{N})\|_{\dot {FB}^{2-\frac{3}{p}}_{p,r}}\leq C\nu,\eqno(1.5)$$
then problem (1.4) possesses a unique global mild solution in the class $\big[C\big([0,\infty);\dot {FB}^{2-\frac{3}{p}}_{p,r}(\mathbb{R}^3)\big)\big]^4\cap X_\infty^\alpha$.

\vskip 3mm

{\bf Theorem 1.5}\, \,  Let Prandtl number $P=1$, i.e., $\nu=\mu$, and $r\in[1, 2]$. Then for any $(u_0,\theta_0)\in\big[\dot{FB}^{-1}_{1,r}(\mathbb{R}^3)\big]^{4}$ satisfying
$\textrm{div} u_0=0$, there exists corresponding $T>0$ such that problem (1.4) possesses a unique mild solution $(u, \theta)\in
\big[C\big([0,T],\dot{FB}^{-1}_{1,r}(\mathbb{R}^3)\big)\big]^{4}\cap Y_T^\alpha$, where
$$Y_T^\alpha:=[\tilde{L}^{\frac{2}{1+\alpha}}(0,T;\dot{FB}^{\alpha}_{1,r}(\mathbb{R}^3))]^4
\cap [\tilde{L}^{\frac{2}{1-\alpha}}(0,T;\dot{FB}^{-\alpha}_{1,r}(\mathbb{R}^3))]^4$$
(for an arbitrarily chosen but fixed number $\alpha\in(0,1)$).
Moreover, there exists a constant $C>0$ independent of $\nu$,
$\Omega$ and $\mathcal{N}$ such that if
$$\|(u_0, \sqrt{g}\theta_0/\mathcal{N})\|_{\dot {FB}^{-1}_{1,r}}\leq C\nu,\eqno(1.6)$$
then problem (1.4) possesses a unique global mild solution in the class $\big[C\big([0,\infty);\dot {FB}^{-1}_{1,r}(\mathbb{R}^3)\big)\big]^4\cap Y_\infty^\alpha$.

\vskip 3mm

{\bf Theorem 1.6 }\, \, Let Prandtl number $P=1$, i.e., $\nu=\mu$ and $2<r\leq \infty$. Then the problem (1.4) is ill-posed in $\big[\dot{FB}^{-1}_{1,r}(\mathbb{R}^3)\big]^4$
in the sense that the solution map $(u_0,\theta_0)\mapsto(u,\theta)$ from $\big[\dot{FB}^{-1}_{1,r}(\mathbb{R}^3)\big]^4$ to
$\big[C\big([0,T],\dot{FB}^{-1}_{1,r}(\mathbb{R}^3)\big)\big]^4$, if exists, is not continuous at $(u_0,\theta_0)=(0,0)$.

\medskip

Theorems 1.4 and 1.5 show that the Cauchy problem (1.4) with Prandtl number $P=1$ is locally well-posed in
the Fourier-Besov spaces $[\dot{FB}^{2-\frac{3}{p}}_{p,r}(\mathbb{R}^3)]^4$ for $1<p\leq \infty,1\leq r< \infty$ and $p=1$, $1\leq r\leq 2$, and
globally well-posed in these spaces when the initial data $(u_0,\theta_0)$ satisfy smallness conditions (1.5) and (1.6). On the other hand, Theorem
1.6 tells us that this problem is ill-posed in $[\dot{FB}^{-1}_{1,r}(\mathbb{R}^3)]^4$ for $2<r\leq\infty$. Hence, we have established the sharp
well-posedness and ill-posedness in Fourier-Besov spaces $[\dot{FB}^{2-\frac{3}{p}}_{p,r}(\mathbb{R}^3)]^4$ for problem (1.4) under
the special case Prantl number $P=1$. The condition $P=1$ is imposed for similar technical reasons as in \cite{Koba}; to the present we are not clear whether
this condition can be removed or not.

The proofs of Theorems 1.4 and 1.5 use the standard Picard iteration argument. The main point is that in order to get the sharp well-posedness results
as stated in these theorems we must be very careful in the construction of the iteration scheme. Indeed, as well known, to get well-posedness
of Cauchy problems of evolution equations in function spaces possessing sufficiently high regularity, there are many different choices of the
iteration scheme. However, most of those schemes do not work in function spaces of low regularity. Concerning the problem (1.4), we shall use the
same iteration scheme as that used in the literature \cite{Koba}. Thanks to the hypothesis $P=1$ or $\mu=\nu$, the semigroup $T_{\Omega,N}$ related
to that iteration scheme satisfies certain very nice estimates similar to those established for the semigroup related to the equations (1.3)
in the literature \cite{Iwabuchi 2014}; see Section 2 for details. The proof of Theorem 1.6 is much harder. We shall use some arguments similar to
those in \cite{Bejenaru} and \cite{Iwabuchi 2014} to prove this theorem. Note that since in the present case we consider the stratification effects of the flow and one more unknown function
$\theta$ than in \cite{Iwabuchi 2014}, the analysis is more involved; see Section 3 for details.

The rest of this paper is organized as follows. In Section 2 we first transform the initial value problem (1.4) into an equivalent integral
equation (which is essentially the same as to construct an iteration scheme), next establish some
linear estimates and product laws, and finally we present the proofs of Theorems 1.4 and 1.5. The last section is devoted to giving the proof of
Theorem 1.6.

Throughout this paper, we shall use $C$ and $c$ to denote universal constants whose value
may change from line to line. Both $\mathcal{F}g$ and $\hat{g}$ stand for Fourier transform
of $g$ with respect to space variable, while $\mathcal{F}^{-1}$ stands for the inverse Fourier transform.
Besides, since we only
consider the case $P=1$, we always assume that $\mu=\nu$.

\vskip 4mm

\section{Proofs of Theorems 1.4 and 1.5}

In this section we present the proofs of Theorems 1.4 and 1.5. To this end, we first transform the Cauchy problem
(1.4) into an equivalent integral equation, and next use the Littlewood-Paley analysis technique to establish some linear estimates and product laws.

\subsection{Rewriting (1.4) into an integral equation}
By setting $N:=\mathcal{N}\sqrt{g}$, $v:=(v^1,v^2,v^3,v^4):=
(u^1,u^2,u^3,\sqrt{g}\theta/\mathcal{N})$,  $v_0:=(v_0^1,v_0^2,v_0^3,v_0^4):=(u_0^1,u_0^2,u_0^3,\sqrt{g}\theta_0/\mathcal{N})$ and
$\widetilde{\nabla}:=(\partial_1,\partial_2,\partial_3,0)$, (1.4) can be rewritten into the following problem
$$
\left\{
  \begin{array}{ll}
   \partial_t v+\mathcal{A} v+\mathcal{B} v+\widetilde{\nabla} p=-(v\cdot\widetilde{\nabla})v
   \ \ \ \ & \textrm{in}~~\mathbb{R}^3\times (0,\infty), \\
   \widetilde{\nabla}\cdot v=0  & \textrm{in}~~\mathbb{R}^3\times (0,\infty), \\
   v|_{t=0}=v_0  & \textrm{in}~~\mathbb{R}^3,
  \end{array}
\right.\eqno(2.1)
$$
where
$$\mathcal{A}:=\begin{pmatrix}
  -\nu \Delta & 0 & 0 & 0 \\
  0 &  -\nu \Delta & 0 & 0 \\
  0 & 0 &  -\nu \Delta & 0\\
  0 & 0 & 0 &  -\mu \Delta \\
\end{pmatrix}~~~~ \mbox{and} ~~~~
\mathcal{B}:=\begin{pmatrix}
  0~~ & -\Omega~~ & 0~~ & 0 \\
  \Omega~~ & 0~~ & 0~~ & 0 \\
  0~~ & 0~~ & 0~~ & -N \\
  0~~ & 0~~ & N~~ & 0 \\
\end{pmatrix}.\eqno(2.2)
$$
Lemma 3.3 in \cite{Koba}, together with the fact $e^{(\mathcal{A}+\mathcal{B})t}=e^{\mathcal{A}t}e^{\mathcal{B}t}$ for $\nu=\mu$, gives the explicitly expression of Stokes-Coriolis-Stratification semigroup $T_{\Omega,N}(t)(t\geq 0)$ corresponding to the linear problem of (2.1) via Fourier transform
$$
  T_{\Omega,N}(t)f:=\mathscr{F}^{-1}\Big[\cos\Big(\frac{|\xi|'}{|\xi|} t\Big)e^{-\nu|\xi|^2t}M_1\hat{f}+\sin\Big(\frac{|\xi|'}{|\xi|} t\Big)
  e^{-\nu|\xi|^2t}M_2\hat{f}+e^{-\nu|\xi|^2t}M_3\hat{f}\Big],
\eqno(2.3)
$$
where
$$
  |\xi|:=\sqrt{\xi_1^2+\xi_2^2+\xi_3^2},~~~~~|\xi|':=|\xi|'_{\Omega, N}:=\sqrt{N^2\xi_1^2+N^2\xi_2^2+\Omega^2\xi_3^2}\eqno(2.4)
$$
for $\xi=(\xi_1,\xi_2,\xi_3)\in\mathbb{R}^3$,
$$
  M_1=\begin{pmatrix}
  \frac{\Omega^2\xi_3^2}{|\xi|'^2}~ & ~0~ & ~-\frac{N^2\xi_1\xi_3}{|\xi|'^2}~ & ~\frac{\Omega N \xi_2\xi_3}{|\xi|'^2}\\
  0 & \frac{\Omega^2\xi_3^2}{|\xi|'^2}~ & ~-\frac{N^2\xi_2\xi_3}{|\xi|'^2}~ & ~-\frac{\Omega N \xi_1\xi_3}{|\xi|'^2} \\
  -\frac{\Omega^2\xi_1\xi_3}{|\xi|'^2}~ & ~-\frac{\Omega^2\xi_2\xi_3}{|\xi|'^2}~ & ~\frac{N^2(\xi_1^2+\xi_2^2)}{|\xi|'^2}~ & ~0 \\
  \frac{\Omega N\xi_2\xi_3}{|\xi|'^2}~ & ~-\frac{\Omega N\xi_1\xi_3}{|\xi|'^2}~ & ~0~ & ~\frac{N^2(\xi_1^2+\xi_2^2)}{|\xi|'^2} \\
  \end{pmatrix},
\eqno(2.5)$$
$$
  M_2=\begin{pmatrix}
    0~ & ~-\frac{\Omega\xi_3^2}{|\xi||\xi|'}~ & ~\frac{\Omega\xi_2\xi_3}{|\xi||\xi|'}~  & ~\frac{N\xi_1\xi_3}{|\xi||\xi|'} \\
    \frac{\Omega\xi^2_3}{|\xi||\xi|'}~ & ~0~ & ~-\frac{\Omega\xi_1\xi_3}{|\xi||\xi|'}~ & ~\frac{N\xi_2\xi_3}{|\xi||\xi|'} \\
    -\frac{\Omega\xi_2\xi_3}{|\xi||\xi|'}~ & ~\frac{\Omega\xi_1\xi_3}{|\xi||\xi|'}~ & ~0~ & ~-\frac{N(\xi_1^2+\xi^2_3)}{|\xi||\xi|'} \\
    -\frac{N\xi_1\xi_3}{|\xi||\xi|'}~ & ~-\frac{N\xi_2\xi_3}{|\xi||\xi|'}~ & ~\frac{N(\xi_1^2+\xi^2_3)}{|\xi||\xi|'}~ & ~0 \\
  \end{pmatrix}
\eqno(2.6)
$$
and
$$
  M_3=\begin{pmatrix}
  \frac{N^2\xi^2_2}{|\xi|'^2}~~ & -\frac{N^2\xi_1\xi_2}{|\xi|'^2}~~ & 0~~ & -\frac{N\Omega\xi_2\xi_3}{|\xi|'^2} \\
  -\frac{N^2\xi_1\xi_2}{|\xi|'^2}~~ & \frac{N^2\xi^2_1}{|\xi|'^2}~~ & 0~~ & \frac{N\Omega\xi_1\xi_3}{|\xi|'^2} \\
  0~~ & 0~~ & 0~~ & 0 \\
  -\frac{N\Omega\xi_2\xi_3}{|\xi|'^2}~~ & \frac{N\Omega\xi_1\xi_3}{|\xi|'^2}~~ & 0~~ & \frac{\Omega^2\xi^2_3}{|\xi|'^2} \\
  \end{pmatrix}.
\eqno(2.7)
$$
Note that, denoting by $M^l_{jk}(\xi)$ the $(j,k)$-th component of the matrix
$M_l(\xi)$, it is obvious that non-vanishing $M^l_{jk}(\xi)$ satisfies
$$
  |M^l_{jk}(\xi)|\leq 2 \quad  \mbox{for}\;\; \xi\in\mathbb{R}^3, \quad j,k=1,2,3,4,\;\; l=1,2,3.
$$
Define the partial Helmholtz projection operator $\widetilde{\mathbb{P}}=
(\widetilde{\mathbb{P}}_{jk})_{4\times 4}$ by
$$
\widetilde{\mathbb{P}}_{jk}:=\left\{
  \begin{array}{ll}
   \delta_{jk}+R_jR_k,~~~~ 1\leq j,k \leq3, \\
   \delta_{jk},~~~~~~~~~~~~~~~\textrm{otherwise},
  \end{array}
\right.
\eqno(2.8)
$$
where $\delta_{jk}$ is the Kronecker's delta notation
and $R_j$ ($j=1,2,3$) are the Riesz transforms on $\mathbb{R}^3$. Then, by using the Duhamel principle, we easily see that problem (2.1) is equivalent to
the following integral equation
$$
  v(t)=T_{\Omega, N}(t)v_0-\int^t_0T_{\Omega, N}(t-\tau)\widetilde{\mathbb{P}}\widetilde{\nabla}\cdot[v(\tau)\otimes v(\tau)]d\tau.
\eqno(2.9)
$$

\subsection{Linear estimates and product laws}

Next we establish some basic estimates which will play a crucial role in the proofs of Theorems 1.4 and 1.5. We first consider linear estimates for
the semigroup $\{T_{\Omega,N}(t)\}_{t\geq 0}$.
\medskip

{\bf Lemma 2.1}   Let $T>0$, $s\in\mathbb{R}$, $p,r\in[1,\infty]$ and $\alpha\in[0,1]$. There exists a constant $C>0$ such that
$$
\|T_{\Omega,N}(\cdot)u_0\|_{\tilde{L}^{\frac{2}{1\pm\alpha}}(0,T;\dot{FB}^{s+1\pm\alpha}_{p,r})}
\leq C\nu^{-\frac{(1\pm\alpha)}{2}}\|u_0\|_{\dot{FB}^s_{p,r}}
$$
for $u_0\in\dot{FB}^s_{p,r}(\mathbb{R}^3)$.
\medskip

{\bf Proof.} Since $\textrm{supp}\hat{\psi}_j\subset
\{\xi\in\mathbb{R}^3:2^{j-1}\leq|\xi|\leq2^{j+1}\}$, one has
\begin{eqnarray}
\|\mathscr{F}[\Delta_j T_{\Omega,N}(\cdot)u_0]\|_{L^p}&\leq& C\bigg\{\int\limits_{2^{j-1}\leq|\xi|
\leq 2^{j+1}}\!\!\!e^{-\nu|\xi|^2tp}|\hat{\psi}_j(\xi)\hat{u}_0(\xi)|^pd\xi\bigg\}^{\frac{1}{p}}\leq
 Ce^{-\nu2^{2j}t}\|\hat{\psi}_j\hat{u}_0\|_{L^p}\nonumber
\end{eqnarray}
for all $t\geq0$, which yields that
$$\|\mathscr{F}[\Delta_j T_{\Omega,N}(\cdot)u_0]\|_{L^{\frac{2}{1\pm\alpha}}(0,T; L^p)}
\leq C\bigg(\frac{1-e^{\nu2^{2j}\frac{2}{1\pm\alpha}T}}{\nu2^{2j}\frac{2}{1\pm\alpha}}\bigg)^{\frac{1\pm\alpha}{2}}
\|\hat{\psi}_j\hat{u}_0\|_{L^p}.$$
Thus, we have
\begin{eqnarray}
\|T_{\Omega,N}(\cdot)u_0\|_{\tilde{L}^{\frac{2}{1\pm\alpha}}(0,T;\dot{FB}^{s+1\pm\alpha}_{p,r})}
&\leq&C\bigg[\sum_{j\in \mathbb{Z}}\bigg(\frac{1\pm\alpha}{2\nu}\bigg)^{\frac{(1\pm\alpha)r}{2}}
(2^{js}\|\hat{\psi}_j\hat{u}_0\|_{L^p})^r\bigg]^{\frac{1}{r}}\nonumber\\
&\leq& C\nu^{-\frac{(1\pm\alpha)}{2}}\|u_0\|_{\dot{FB}^s_{p,r}}.\nonumber
\end{eqnarray}
\rightline{$\Box$}

{\bf Lemma 2.2} Let $T>0$, $s\in\mathbb{R}$, $1\leq p,r\leq \infty$ and $\alpha\in[0,1]$. There exists a constant $C>0$ such that
$$\bigg\|\int_0^tT_{\Omega, N}(t-\tau)f(\tau)d\tau\bigg\|_{\tilde{L}^{\frac{2}{1\pm\alpha}}(0,T;\dot{FB}^{s+1\pm\alpha}_{p,r})}
\leq C\nu^{-(1+\frac{1\pm\alpha}{2}-\frac{1}{\rho})}
\|f\|_{\tilde{L}^{\rho}(0,T;\dot{FB}^{s-2+\frac{2}{\rho}}_{p,r})}$$
for $f\in\tilde{L}^{\rho}(0,T;\dot{FB}^{s-2+\frac{2}{\rho}}_{p,r}(\mathbb{R}^3)) $ with $\rho\in[1,\frac{2}{1\pm\alpha}]$.
\medskip

{\bf Proof.} By the definition of the $\tilde{L}^{\frac{2}{1\pm\alpha}}(0,T;\dot{FB}^{s+1\pm\alpha}_{p,r}(\mathbb{R}^3))$ and by Young's
inequality, one has
\begin{eqnarray}
&&~\bigg\|\int_0^tT_{\Omega, N}(t-\tau)f(\tau)d\tau\bigg\|_{\tilde{L}^{\frac{2}{1\pm\alpha}}
(0,T;\dot{FB}^{s+1\pm\alpha}_{p,r})}\nonumber\\
&\leq&C\Big\|\Big\{2^{j(s+1\pm\alpha)}\int^t_0e^{-\nu(t-\tau)2^{2j}}
\|\hat{f}(\tau)\cdot\hat{\psi_j}\|_{L^p}
d\tau\bigg\|_{L^{\frac{2}{1\pm\alpha}}(0,T)}\Big\}_{r\in\mathbb{Z}}
\Big\|_{l^r(\mathbb{Z})}\nonumber\\
&\leq&C\Big\|\Big\{2^{j(s+1\pm\alpha)}\|e^{-\nu t2^{2j}}\|_{L^m(0,T)}\|\hat{f}(\tau)\cdot\hat{\psi_j}\|_{L^{\rho}(0,T;L^p)}
\Big\}_{r\in\mathbb{Z}}\Big\|_{l^r(\mathbb{Z})},\nonumber
\end{eqnarray}
where $1+\frac{1\pm\alpha}{2}=\frac{1}{\rho}+\frac{1}{m}$. Thus, we obtain
\begin{eqnarray}
&&~\bigg\|\int_0^tT_{\Omega, N}(t-\tau)f(\tau)d\tau\bigg\|_{\tilde{L}^{\frac{2}{1\pm\alpha}}
(0,T;\dot{FB}^{s+1\pm\alpha}_{p,r})}\nonumber\\
&\leq&C\nu^{-(1+\frac{1\pm\alpha}{2}-\frac{1}{\rho})}
\bigg(\sum_{j\in\mathbb{Z}}2^{j(s+1\pm\alpha)r}2^{-2j(1+\frac{1\pm\alpha}{2}-\frac{1}{\rho})r}
\|\hat{f}(\tau)\cdot\hat{\psi_j}\|_{L^{\rho}(0,T;L^p)}^r\bigg)^\frac{1}{r}\nonumber\\
&\leq& C\nu^{-(1+\frac{1\pm\alpha}{2}-\frac{1}{\rho})}
\|f\|_{\tilde{L}^{\rho}(0,T;\dot{FB}^{s-2+\frac{2}{\rho}}_{p,r})}.\nonumber
\end{eqnarray}
\rightline{$\Box$}
\vskip 2mm

We now turn to establish product laws.
\medskip

{\bf Lemma 2.3} Let $T>0$, $\alpha\in(0,1]$, $p\in(1,\infty]$ and $r\in[1,\infty]$. There exists a constant $C>0$ such that
\begin{eqnarray}
\|fg\|_{\tilde{L}^{1}(0,T; \dot {FB}^{3-\frac{3}{p}}_{p,r})}&\leq & C\bigg(\|f\|_{\tilde{L}^{\frac{2}{1+\alpha}}(0,T;
\dot {FB}^{3-\frac{3}{p}+\alpha}_{p,r})}\|g\|_{\tilde{L}^{\frac{2}{1-\alpha}}(0,T; \dot {FB}^{3-\frac{3}{p}-\alpha}_{p,r})}\nonumber\\
&&~~~~~~~~~~~~+\|g\|_{\tilde{L}^{\frac{2}{1+\alpha}}(0,T; \dot {FB}^{3-\frac{3}{p}+\alpha}_{p,r})}
\|f\|_{\tilde{L}^{\frac{2}{1-\alpha}}(0,T; \dot {FB}^{3-\frac{3}{p}-\alpha}_{p,r})}\bigg)
\nonumber
\end{eqnarray}
for all $f,g\in\tilde{L}^{\frac{2}{1\pm\alpha}}(0,T; \dot {FB}^{3-\frac{3}{p}\pm\alpha}_{p,r}(\mathbb{R}^3))$.
\medskip

{\bf Proof.} Bony's decomposition (see \cite{Bahouri, Bony}) for $\Delta_j(fg)$ reads
\begin{eqnarray}
\Delta_j(fg)&=&\sum_{|k-j|\leq 4}\Delta_jS_{k-1}f\Delta_{k}g
+\sum_{|k-j|\leq 4}\Delta_jS_{k-1}f\Delta_{k}g
+\sum_{k\geq j-2}\sum_{|k'-k|\leq1}\Delta_j\Delta_{k}f\Delta_{k'}g\nonumber\\
&:=&I_1+I_2+I_3.\setcounter{equation}{9}
\end{eqnarray}
Then by the triangle inequalities in $l^p(\mathbb{Z})$ and $L^p(\mathbb{R})$, we have
\begin{eqnarray}
\|fg\|_{\tilde{L}^{1}(0,T; \dot {FB}^{3-\frac{3}{p}}_{p,r})}
&\leq&\bigg(\sum_{j\in\mathbb{Z}}2^{j(3-\frac{3}{p})r}
\|\widehat{I_1}\|^r_{L^{1}(0,T;L^p)}\bigg)^{\frac{1}{r}}
+\bigg(\sum_{j\in\mathbb{Z}}2^{j(3-\frac{3}{p})r}
\|\widehat{I_2}\|^r_{L^{1}(0,T;L^p)}\bigg)^{\frac{1}{r}}\nonumber\\
&&~+\bigg(\sum_{j\in\mathbb{Z}}2^{j(3-\frac{3}{p})r}
\|\widehat{I_3}\|^r_{L^{1}(0,T;L^p)}\bigg)^{\frac{1}{r}}\nonumber\\
&:=&J_1+J_2+J_3.\nonumber
\end{eqnarray}
For $J_1$, Young's inequality and H\"{o}lder's inequality ensures that
\begin{eqnarray}
2^{j(3-\frac{3}{p})}
\|\widehat{I_1}\|_{L^{1}(0,T;L^p)}
&\leq&2^{j(3-\frac{3}{p})}\bigg\|\sum_{|k-j|\leq 4}\big\|\hat{\psi}_j[(\sum_{k'\leq k-2}\hat{\psi}_{k'}\hat{f})\ast(\hat{\psi}_{k}\hat{g})]\big\|_{L^p}\bigg\|_{L^1(0,T)}\nonumber\\
&\leq&C2^{j(3-\frac{3}{p})}\sum_{|k-j|\leq 4}\bigg(\sum_{k'\leq k-2}2^{k'(3-\frac{3}{p})}\big\|\hat{\psi}_{k'}\hat{f}\big\|_{L^{\frac{2}{1-\alpha}}(0,T;L^p)}\bigg)
\big\|\hat{\psi}_{k}\hat{g}\big\|_{L^{\frac{2}{1+\alpha}}(0,T;L^p)}\nonumber\\
&\leq&C2^{j(3-\frac{3}{p})}\sum_{|k-j|\leq 4}\bigg(\sum_{k'\leq k-2}2^{k'(3-\frac{3}{p}-\alpha)r}
\big\|\hat{\psi}_{k'}\hat{f}\big\|^r_{L^{\frac{2}{1-\alpha}}(0,T;L^p)}\bigg)^{\frac{1}{r}}\nonumber\\
&&~~~~~~~~~~~~~~~~~~~~~~~~~~\times\bigg(\sum_{k'\leq k-2}2^{k'\alpha r'}\bigg)^{\frac{1}{r'}}
\big\|\hat{\psi}_{k}\hat{g}\big\|_{L^{\frac{2}{1+\alpha}}(0,T;L^p)}\nonumber\\
&\leq&C\sum_{|k-j|\leq 4}2^{(j-k)(3-\frac{3}{p})}2^{k(3-\frac{3}{p}+\alpha)}
\big\|\hat{\psi}_{k}\hat{g}\big\|_{L^{\frac{2}{1+\alpha}}(0,T;L^p)}
\|f\|_{\tilde{L}^{\frac{2}{1-\alpha}}(0,T; \dot {FB}^{3-\frac{3}{p}-\alpha}_{p,r})}\nonumber
\end{eqnarray}
since $\alpha>0$, where $\frac{1}{r}+\frac{1}{r'}=1$. Hence, by the Young's inequality, we get
$$J_1\leq C\|f\|_{\tilde{L}^{\frac{2}{1-\alpha}}(0,T; \dot {FB}^{3-\frac{3}{p}-\alpha}_{p,r})}\|g\|_{\tilde{L}^{\frac{2}{1+\alpha}}(0,T; \dot {FB}^{3-\frac{3}{p}+\alpha}_{p,r})}.$$
Similarly, we have
$$J_2\leq C\|f\|_{\tilde{L}^{\frac{2}{1+\alpha}}(0,T; \dot {FB}^{3-\frac{3}{p}+\alpha}_{p,r})}\|g\|_{\tilde{L}^{\frac{2}{1-\alpha}}(0,T; \dot {FB}^{3-\frac{3}{p}-\alpha}_{p,r})}.$$
For $J_3$, Young's inequality together with H\"{o}lder's inequality gives
\begin{eqnarray}
2^{j(3-\frac{3}{p})}
\|\widehat{I_3}\|_{L^{1}(0,T;L^p)}
&\leq&2^{j(3-\frac{3}{p})}\bigg\|\sum_{k\geq j-2}\big\|\hat{\psi}_j[(\hat{\psi}_{k}\hat{f})\ast(\sum_{|k'- k|\leq1}\hat{\psi}_{k'}\hat{g})]\big\|_{L^p}\bigg\|_{L^1(0,T)}\nonumber\\
&\leq&C2^{j(3-\frac{3}{p})}\sum_{k\geq j-2}2^{3k(1-\frac{1}{p})}\big\|\hat{\psi}_{k}\hat{f}\big\|_{L^{\frac{2}{1+\alpha}}(0,T;L^p)}
\bigg(\sum_{|k'-k|\leq1}2^{-k'(3-\frac{3}{p}-\alpha)r'}\bigg)^{\frac{1}{r'}}
\nonumber\\
&&~~~~~~~~~~~~~~~~~~~~\times\bigg(\sum_{|k'-k|\leq1}2^{k'(3-\frac{3}{p}-\alpha)r}
\big\|\hat{\psi}_{k'}\hat{g}\big\|^r_{L^{\frac{2}{1-\alpha}}(0,T;L^p)}\bigg)^{\frac{1}{r}}
\nonumber\\
&\leq&C\sum_{k\geq j-2}2^{(j-k)(3-\frac{3}{p})}2^{k(3-\frac{3}{p}+\alpha)}
\big\|\hat{\psi}_{k}\hat{f}\big\|_{L^{\frac{2}{1+\alpha}}(0,T;L^p)}
\|g\|_{\tilde{L}^{\frac{2}{1-\alpha}}(0,T; \dot {FB}^{3-\frac{3}{p}-\alpha}_{p,r})},\nonumber
\end{eqnarray}
where $\frac{1}{r}+\frac{1}{r'}=1$. Hence, by the Young's inequality, one has
\begin{eqnarray}J_3&\leq& C\bigg(\sum_{2\geq k}2^{k(3-\frac{3}{p})}\bigg) \|f\|_{\tilde{L}^{\frac{2}{1-\alpha}}(0,T; \dot {FB}^{3-\frac{3}{p}-\alpha}_{p,r})}\|g\|_{\tilde{L}^{\frac{2}{1+\alpha}}(0,T; \dot {FB}^{3-\frac{3}{p}+\alpha}_{p,r})}\nonumber\\
&\leq& C\|f\|_{\tilde{L}^{\frac{2}{1-\alpha}}(0,T; \dot {FB}^{3-\frac{3}{p}-\alpha}_{p,r})}\|g\|_{\tilde{L}^{\frac{2}{1+\alpha}}(0,T; \dot {FB}^{3-\frac{3}{p}+\alpha}_{p,r})}\nonumber
\end{eqnarray}
since $p>1$.
Summing up, we arrive at
\begin{eqnarray}
\|fg\|_{\tilde{L}^{1}(0,T; \dot {FB}^{3-\frac{3}{p}}_{p,r})}&\leq & C\bigg(\|f\|_{\tilde{L}^{\frac{2}{1+\alpha}}(0,T; \dot {FB}^{3-\frac{3}{p}+\alpha}_{p,r})}\|g\|_{\tilde{L}^{\frac{2}{1-\alpha}}(0,T; \dot {FB}^{3-\frac{3}{p}-\alpha}_{p,r})}\nonumber\\
&&~~~~~~~~~~~~+\|g\|_{\tilde{L}^{\frac{2}{1+\alpha}}(0,T; \dot {FB}^{3-\frac{3}{p}+\alpha}_{p,r})}\|f\|_{\tilde{L}^{\frac{2}{1-\alpha}}(0,T; \dot {FB}^{3-\frac{3}{p}-\alpha}_{p,r})}\bigg)\nonumber
\end{eqnarray}
\rightline{$\Box$}

The above lemma excludes the end point case $p=1$, which is considered in the following one.
\medskip

{\bf Lemma 2.4} Let $T>0$, $\alpha\in(0,1]$ and $r\in[1,2]$. There exists a constant $C>0$ such that
\begin{eqnarray}
\|fg\|_{\tilde{L}^{1}(0,T; \dot {FB}^{0}_{1,r})}&\leq & C\bigg(\|f\|_{\tilde{L}^{\frac{2}{1+\alpha}}(0,T;
\dot {FB}^{\alpha}_{1,r})}\|g\|_{\tilde{L}^{\frac{2}{1-\alpha}}(0,T; \dot {FB}^{-\alpha}_{1,r})}\nonumber\\
&&~~~~~~~~~~~~+\|g\|_{\tilde{L}^{\frac{2}{1+\alpha}}(0,T; \dot {FB}^{\alpha}_{1,r})}
\|f\|_{\tilde{L}^{\frac{2}{1-\alpha}}(0,T; \dot {FB}^{-\alpha}_{1,r})}\bigg)\nonumber
\end{eqnarray}
for all $f,g\in\tilde{L}^{\frac{2}{1\pm\alpha}}(0,T; \dot {FB}^{\pm\alpha}_{1,r}(\mathbb{R}^3))$.
\medskip

{\bf Proof} Similar to the proof of Lemma 2.3, by means of Bony's decomposition we have
\begin{eqnarray}
\|fg\|_{\tilde{L}^{1}(0,T; \dot {FB}^{0}_{1,r})}
&\leq&\bigg(\sum_{j\in\mathbb{Z}}
\|\widehat{I_1}\|^r_{L^{1}(0,T;L^1)}\bigg)^{\frac{1}{r}}
+\bigg(\sum_{j\in\mathbb{Z}}
\|\widehat{I_2}\|^r_{L^{1}(0,T;L^1)}\bigg)^{\frac{1}{r}}\nonumber\\
&&~+\bigg(\sum_{j\in\mathbb{Z}}
\|\widehat{I_3}\|^r_{L^{1}(0,T;L^1)}\bigg)^{\frac{1}{r}}\nonumber\\
&:=&J_1+J_2+J_3,\nonumber
\end{eqnarray}
where $I_i(i=1,2,3)$ are defined in (2.9).
It is easy to check that the estimates for $J_1$ and $J_2$ in the proof of Lemma 2.3 also hold for the case that $p=1$ with $\alpha\in(0,1]$. Thus, we
only need to consider $J_3$.

Applying H\"{o}lder's inequality and Young's inequality gives
\begin{eqnarray}
J_3
&\leq&
\big\|\sum_{k\in\mathbb{Z}}\big\|[(\hat{\psi}_{k}\hat{f})\ast(\sum_{|k'- k|\leq1}\hat{\psi}_{k'}\hat{g})]\big\|_{L^1}\big\|_{L^1(0,T)}\nonumber\\
&\leq&\bigg(\sum_{k\in\mathbb{Z}}2^{-k\alpha r}
\big\|\hat{\psi}_{k}\hat{f}\big\|^r_{L^{\frac{2}{1-\alpha}}(0,T;L^1)}\bigg)^{\frac{1}{r}}
\bigg\|\bigg(\sum_{|k'-k|\leq1}2^{(k-k')\alpha}2^{k'\alpha}
\big\|\hat{\psi}_{k'}\hat{g}\big\|_{L^{\frac{2}{1+\alpha}}(0,T; L^1)}\bigg)\bigg\|_{l^{r'}}
\nonumber\\
&\leq&\|f\|_{\tilde{L}^{\frac{2}{1-\alpha}}(0,T; \dot {FB}^{-\alpha}_{1,r})}\bigg(\sum_{|k-k'|\leq1}
2^{k'\alpha r}\big\|\hat{\psi}_{k'}\hat{g}\big\|^r_{L^{\frac{2}{1+\alpha}}(0,T;L^1)}\bigg)^{\frac{1}{r}}
\bigg(\sum_{|k|\leq1}2^{k\alpha m}\bigg)^{\frac{1}{m}}
\nonumber\\
&\leq&C \|f\|_{\tilde{L}^{\frac{2}{1-\alpha}}(0,T; \dot {FB}^{-\alpha}_{1,r})}
\|g\|_{\tilde{L}^{\frac{2}{1+\alpha}}(0,T; \dot {FB}^{\alpha}_{1,r})},
\nonumber
\end{eqnarray}
where the first inequality has used $l^1\hookrightarrow l^r$ and $\sum_{j\in\mathbb{Z}}\hat{\psi}_{j}=1$, $\frac{1}{r}+\frac{1}{r'}=1$ and $m\in[1,\infty]$ satisfying $\frac{1}{r}+\frac{1}{m}=1+\frac{1}{r'}$ for $r\in[1,2]$.
\rightline{$\Box$}
\vskip 2mm

\subsection{Proofs of Theorems 1.4 and 1.5}

The proofs of Theorems 1.4 and 1.5 follow from the following standard Banach fixed point lemma combined with the estimates
established in the previous section.
\medskip

{\bf  Lemma 2.5} (Cannone and Karch \cite{Cannone 2004})\ \
Let $(\mathcal{X},\|\cdot\|_{\mathcal{X}})$ be a Banach space and $B:\mathcal{X}\times\mathcal{X}\to\mathcal{X}$ a bounded bilinear form satisfying
$\|B(x_1,x_2)\|_{\mathcal{X}}\leq \eta\|x_1\|_{\mathcal{X}}\|x_2\|_{\mathcal{X}}$ for all $x_1, x_2\in\mathcal{X}$ and some constant $\eta>0$.
Then, if $0<\varepsilon<\frac{1}{4\eta}$ and if $y\in\mathcal{X}$ such that $\|y\|_{\mathcal{X}}<\varepsilon$, the equation $x=y+B(x,x)$ has a
solution in $\mathcal{X}$ such that $\|x\|_{\mathcal{X}}\leq2\varepsilon$.
This solution is the only one in the ball $\bar{B}(0,2\varepsilon)$.
Moreover, the solution depends continuously on $y$ in the following sense: if $\|\tilde{y}\|_{\mathcal{X}}\leq\varepsilon$, $\tilde{x}=
\tilde{y}+B(\tilde{x},\tilde{x})$ and $\|\tilde{x}\|_{\mathcal{X}}\leq2\varepsilon$ then
$$\|x-\tilde{x}\|_{\mathcal{X}}\leq \frac{1}{1-4\eta\varepsilon}\|y-\tilde{y}\|_{\mathcal{X}}.$$
\rightline{$\Box$}

{\bf  Proof of Theorem 1.4}\ \ \ Define
$$B(v,w)(t):=\int^t_0T_{\Omega, N}(t-\tau)\widetilde{\mathbb{P}}\widetilde{\nabla}\cdot[v(\tau)\otimes w(\tau)]d\tau.\eqno(2.10)$$
Let $\alpha\in(0,1)$ be any given and fixed. For $T>0$ to be specified later, let $X^{\alpha}_T$ be the function
space introduced in Theorem 1.4. It is clear that $X^{\alpha}_T$ is a Banach space endowed with the norm
$$\|v\|_{X^{\alpha}_T}:=\|v\|_{\tilde{L}^{\frac{2}{1-\alpha}}(0,T;\dot{FB}^{3-\frac{3}{p}-\alpha}_{p,r})}
+\|v\|_{\tilde{L}^{\frac{2}{1+\alpha}}(0,T;\dot{FB}^{3-\frac{3}{p}+\alpha}_{p,r})}.$$
By applying Lemmas 2.2 and 2.3 with $s=2-\frac{3}{p}$ and $\rho=1$, we arrive at
\begin{eqnarray}
\|B(v,w)\|_{X^{\alpha}_T}&=&\bigg\|\int^t_0T_{\Omega, N}(t-\tau)\widetilde{\mathbb{P}}\widetilde{\nabla}\cdot
[v(\tau)\otimes w(\tau)]d\tau\bigg\|_{X^{\alpha}_T}\nonumber\\
&\leq&C_0(\nu^{-\frac{1-\alpha}{2}}+\nu^{-\frac{1+\alpha}{2}})
\|\widetilde{\nabla}\cdot[v(\tau)\otimes w(\tau)]
\|_{\tilde{L}^{1}(0,T;\dot{FB}^{2-\frac{3}{p}}_{p,r})}\nonumber\\
&\leq&C_0\max\{\nu^{-\frac{1+\alpha}{2}},\nu^{-\frac{1-\alpha}{2}}\}\|v\|_{X^{\alpha}_T}\|w\|_{X^{\alpha}_T}\nonumber
\end{eqnarray}
for $v,w\in{X^{\alpha}_T}$ and some constant $C_0>0$. Now let $v_0\in \dot{FB}^{2-\frac{3}{p}}_{p,r}(\mathbb{R}^3)$ be given, then there
exists $T>0$ such that
$$\|T_{\Omega,N}(t)v_0\|_{X^{\alpha}_T}\leq
\frac{\min\{\nu^{\frac{1+\alpha}{2}},\nu^{\frac{1-\alpha}{2}}\}}{8C_0}.\eqno(2.11)$$
With $T>0$ specified in this way, Lemma 2.5 ensures that that there exists a unique solution $v$ of (2.9)
in the ball with center $0$ and radius $\frac{\min\{\nu^{\frac{1+\alpha}{2}},\nu^{\frac{1-\alpha}{2}}\}}{2C_0}$ in the space $X^{\alpha}_T$.
Moreover, applying Lemmas 2.1 $\sim$ 2.3 with $s=2-\frac{3}{p}$ and $\rho=1$ leads to
\begin{eqnarray}
\|v\|_{\tilde{L}^\infty(0,T;\dot{FB}^{2-\frac{3}{p}}_{p,r})}&\leq&C \|v_0\|_{\dot{FB}^{2-\frac{3}{p}}_{p,r}}+C
\|\widetilde{\nabla}\cdot[v(\tau)\otimes v(\tau)]
\|_{\tilde{L}^{1}(0,T;\dot{FB}^{2-\frac{3}{p}}_{p,r})}\nonumber\\
&\leq&C\|v_0\|_{\dot{FB}^{2-\frac{3}{p}}_{p,r}}+C\|v\|^2_{X^{\alpha}_T}<\infty.\nonumber
\end{eqnarray}
By using a standard density
argument, we can further infer that $v\in \big[C([0,T],\dot{FB}^{2-\frac{3}{p}}_{p,r}(\mathbb{R}^3))\big]^4$. This proves the local well-posedness assertion
in Theorem 1.4.

Next we assume that the condition (2.11) is satisfied. It follows from Lemma 2.1 with $s=2-\frac{3}{p}$ that there exists a constant $C_1>0$ such
that
$$\|T_{\Omega,N}(t)v_0\|_{X_T^{\alpha}}\leq C_1\max\{\nu^{-\frac{(1+\alpha)}{2}},\nu^{-\frac{(1-\alpha)}{2}}\}
\|v_0\|_{\dot{FB}^{2-\frac{3}{p}}_{p,r}}$$
for any given $T>0$.
Hence, if
$$\|u_0\|_{\dot{FB}^{2-\frac{3}{p}}_{p,r}}\leq\frac{\nu}{4C_0C_1},$$
then the smallness condition (2.11) with $X^{\alpha}_T$ replaced by $X_\infty^{\alpha}$ is satisfied. Then by deducing similarly as above we see that
problem (2.9) has a unique solution $v\in \big[C\big([0,\infty);\dot {FB}^{2-\frac{3}{p}}_{p,r}(\mathbb{R}^3)\big)\big]^4\cap X_\infty^\alpha$. This proves
the global well-posedness assertion in Theorem 1.4 and finishes the proof of Theorem 1.4.

\rightline{$\Box$}

{\bf  Proof of Theorem 1.5}:~~ The proof of Theorem 1.5 is quite similar to that of Theorem 1.4. What we only need to modify is to replace the
function spaces $X^{\alpha}_T$ and $X_\infty^{\alpha}$ respectively with $Y^{\alpha}_T$ and $Y_\infty^{\alpha}$ introduced in Theorem 1.5.  We omit the details.
\rightline{$\Box$}

\section{Proof of Theorem 1.6}

In this section we give the proof of Theorem 1.6. Similar to \cite{Iwabuchi 2014}, we shall use the following abstract result of Bejenaru and Tao
\cite{Bejenaru} to prove Theorem 1.6.

Consider the abstract equation
$$u=L(f)+N_k(u,\cdots,u),\eqno(3.1)$$
where the initial data $f$ takes values in some data space $D$, the solution $u$ takes values in some solution space $S$, the linear operator
$L: \mathscr{D}(L)\subseteq D\to S$ is densely defined, and the $k$-linear operator $N_k: \mathscr{D}(N_k)\subseteq S\times\cdots\times S\to S$ with
$k\geq2$ is also densely defined. And let $(A_{n_1}(f),\cdots,A_{n_k}(f))\in\mathscr{D}(N_k)$ for all $f\in\mathscr{D}(L)$, $k\geq 2$ and $n_1,n_2,
\cdots,n_k\in\mathbb{N}$, where
$$
\left\{
\begin{array}{l}
  A_1(f):=L(f)\\ [0.2cm]
  A_n(f):=\displaystyle\!\!\sum_{n_1+\cdots+n_k=n}\!\!N_k(A_{n_1}(f),\cdots,A_{n_k}(f))~~~\text{for}~n=2,3,\cdots.
\end{array}
\right.
$$

{\bf Proposition 3.1} (Bejenaru and Tao \cite{Bejenaru}) Suppose that Eq. (3.1) is quantitatively well-posed in the Banach spaces
$D$ and $S$ in the sense that there exists a solution map $f\mapsto u[f]$ from a ball $B_D$ in $D$ to a ball $B_S$ in $S$ which is continuous
with respect to the norm topologies of $D$ and $S$. Suppose that these spaces are endowed with different norms to form different normed vector spaces
$D'$ and $S'$ (not necessarily complete), respectively, which are weaker than $D$ and $S$ in the sense that there exists a constant $C>0$ such
that
$$\|f\|_{D'}\leq C\|f\|_{D} \;\; (\mbox{for}\; f\in D) ~~~~ \mbox{and} ~~~~ \|u\|_{S'}\leq C\|u\|_{S} \;\; (\mbox{for}\; u\in S).$$
Suppose further that the solution map $f\mapsto u[f]$ is continuous from $(B_D, \|\cdot\|_{D'} )$ (i.e. the ball $B_D$ equipped with the $D'$ topology) to $(B_S, \|\cdot\|_{S'})$. Then
for each $n$, the non-linear operator $A_n$ is continuous from $(B_D, \|\cdot\|_{D'})$ to $(S, \|\cdot\|_{S'})$.
\rightline{$\Box$}
\medskip

{\bf  Proof of Theorem 1.6}: ~~ We shall use Proposition 3.1 with the help of Theorem 1.5 to prove Theorem 1.6. Define $D:=
[\dot {FB}^{-1}_{1,2}(\mathbb{R}^3)]^4$ and $S:=\big[C\big([0,\infty);\dot {FB}^{-1}_{1,2}(\mathbb{R}^3)\big)\big]^4\cap [\tilde{L}^{\frac{2}{1+\alpha}}(0,\infty;\dot{FB}^{\alpha}_{1,2}(\mathbb{R}^3))]^4
\cap [\tilde{L}^{\frac{2}{1-\alpha}}(0,\infty;\dot{FB}^{-\alpha}_{1,2}(\mathbb{R}^3))]^4$.  Theorem 1.5 implies that
that for any $0<\alpha<1$ there exists corresponding $\varepsilon>0$ such that there exists a solution map $B_D\ni f\mapsto v[f]\in
S$ which is continuous with respect to the norm topologies of $D$ and $S$, where
$$B_D:=\{f\in D~\mid~\|f\|_D\leq \varepsilon\} ~~~~ \mbox{and} ~~~~  v[f]:=T_{\Omega, N}f-B(v[f],v[f]),$$
where $B(v,w)$ is defined in (2.10).
Let $D':=[\dot {FB}^{-1}_{1,r}(\mathbb{R}^3)]^4$ and $S':=\big[L^\infty(0, \infty; \dot {FB}^{-1}_{1,r})\big]^4$ with $2<r\leq\infty$. It is obvious that the embeddings $D\hookrightarrow D'$ and $S\hookrightarrow S'$ are continuous.

We shall prove by contradiction that the solution map $(B_D,\|\cdot\|_{D'})\ni f\mapsto
v[f]\in(S, \|\cdot\|_{S'})$ is not continuous, no matter how small $\varepsilon$ is chosen.
Hence if we assume that $(B_D, \|\cdot\|_{D'})\ni
f\mapsto v[f]\in(S, \|\cdot\|_{S'})$ is continuous for some $\varepsilon>0$. Then, by Proposition 3.1, the map $(B_D, \|\cdot\|_{D'})\ni f\mapsto
A_2(f)\in(S, \|\cdot\|_{S'})$ is also continuous, where
$$A_2(f):=B\big(T_{\Omega, N}(\cdot)f, T_{\Omega, N}(\cdot)f\big).$$
However, in what follows we shall construct a sequence $\{f^M\}_{M=1}^{\infty}\in B_D$ such that
$$\|f^M\|_{D'}\to 0~~~~\textrm{as}~~M\to \infty,\eqno(3.2)$$
and there exists a constant $c>0$ independent of $M\in\mathbb{N}$ such that
$$\|A_2(f^M)\|_{S'}=\big\|B\big(T_{\Omega, N}(\cdot)f^M, T_{\Omega, N}(\cdot)f^M\big)\big\|_{S'}\geq c\eqno(3.3)$$
for all sufficiently large $M$, and a contradiction follows.

We now construct our counterexample $\{f^M\}_{M=1}^\infty$. Define
$$\chi(\xi)=\left\{
  \begin{array}{ll}
    ~1,~~~~~~~~ & \hbox{\textrm{if}}~|\xi_k|\leq1, k=1,2,3, \\
    ~0,~~~~~~~~ & \hbox{\textrm{otherwise}},
  \end{array}
\right.$$
and $\chi_j^{\pm}(\xi)=\chi(\xi\mp2^je_2)$ for $j\in\mathbb{Z}$, where $e_2=(0,1,0)$.
Then, we set the initial data $\{f^M\}_{M=1}^\infty$ via Fourier transform as follows:
$$\widehat{f^M}(\xi):=\frac{i}{M^{\frac{1}{2}}}\sum_{j=M}^{2M}2^j
\big(\chi_j^+(\xi)+\chi_j^-(\xi)\big)\frac{1}{|\xi|}
\begin{pmatrix}
\xi_2 \\
-\xi_1 \\
0 \\
0 \\
\end{pmatrix}.
$$
By the definition of norm of $\dot{FB}_{1,r}^{-1}$, there exists $C>0$ such that
$$\|f^M\|_{\dot{FB}_{1,r}^{-1}}\leq CM^{-\frac{1}{2}+\frac{1}{r}}, \ \ \ \ \text{for}\ \ M\in\mathbb{N}\ \ \text{and}\ \ 2\leq r\leq \infty.$$
Thus $\varepsilon C^{-1}f^M\in B_D$ and (3.2)
is satisfied for all $2<r\leq \infty$.
Hereinafter, we omit the absolute constant $\varepsilon C^{-1}$ for simplicity and prove the property (3.3).

Let $E$ be a measurable set in $\mathbb{R}^3$ such that the Lebesgue measure of $E$ is positive, there exists a constant $c>0$ such that
$$1-\frac{\xi_1^2}{|\xi|^2}\geq c,\ \ \ \ \ \text{for}\ \ \ \xi=(\xi_1,\xi_2,\xi_3)\in E,\eqno(3.4)$$
and
$$E\subset\bigg\{\xi\in\mathbb{R}^3~\big{|}~\frac{1}{1000}\leq\xi_1\leq1,~ |\xi|\leq1\bigg\}.\eqno(3.5)$$
Since $E$ is bounded, there exists $j_0\in\mathbb{N}$ such that $\sum_{j=-j_0}^{j_0}\hat{\psi}_j(\xi)=1$
for $\xi\in E$. Thus there exists a constant $C_E>0$ such that
\begin{eqnarray}
\|A_2(f^M)(t)\|_{\dot{FB}^{-1}_{1,r}}
&\geq& C_E
\|\mathscr{F}[A_2(f^M)](t)\|_{L^1(E)},\ \ \ \ \ \text{for all}\ \ t\geq 0.\setcounter{equation}{6}
\end{eqnarray}

In what follows, we just need to estimate $\|\mathscr{F}[A_2(f^M)](t)\|_{L^1(E)}$. Now, we assume that $N>0$. If $N=0$, combining with $\theta_0=(f^M)_4=0$, we have $\theta\equiv0$. Then it becomes Cauchy problem of (1.3) which have discussed in \cite{Iwabuchi 2014}. Furthermore, we assume that $N\geq |\Omega|\geq0$. The argument of the case $|\Omega|\geq N>0$ is similar.

For the sake of convenience in writing, we define $\hat{T}_i$, $i=1,2,3$ by
$$\hat{T}_1(\xi,t):=\cos\bigg(\frac{|\xi|'}{|\xi|}t\bigg)e^{-\nu|\xi|^2t},~~~
\hat{T}_2(\xi,t):=\sin\bigg(\frac{|\xi|'}{|\xi|}t\bigg)e^{-\nu|\xi|^2t},\eqno(3.7)$$
and
$$\hat{T}_3(\xi,t):=e^{-\nu|\xi|^2t},\eqno(3.8)$$
where $|\xi|$ and $|\xi|'$ are defined as in (2.4).
Moreover, we have the following key observations:

$\bullet$ $\widetilde{\mathbb{P}}: [\mathscr{S}'(\mathbb{R}^3)]^4\mapsto\mathscr{S}'_\sigma(\mathbb{R}^3)
\times\mathscr{S}'(\mathbb{R}^3)$;

$\bullet$ $(M_1+M_3)(\hat{v})=\hat{v}~~
\text{for}~~v\in\mathscr{S}'_\sigma(\mathbb{R}^3)\times\mathscr{S}'(\mathbb{R}^3)$, since
$$M_1+M_3=\begin{pmatrix}
  1-\frac{N^2\xi_1^2}{|\xi|'^2}~~ & -\frac{N^2\xi_1\xi_2}{|\xi|'^2}~~ & -\frac{N^2\xi_1\xi_3}{|\xi|'^2}~~ & 0 \\
  -\frac{N^2\xi_1\xi_2}{|\xi|'^2}~~ & 1-\frac{N^2\xi_2^2}{|\xi|'^2}~~ & -\frac{N^2\xi_2\xi_3}{|\xi|'^2}~~ & 0 \\
  -\frac{\Omega^2\xi_1\xi_3}{|\xi|'^2}~~ &  -\frac{\Omega^2\xi_2\xi_3}{|\xi|'^2}~~ & 1-\frac{\Omega^2\xi_3^2}{|\xi|'^2}~~ & 0 \\
  0~~ & 0~~ & 0~~ & 1 \\
\end{pmatrix},
$$
where $\mathscr{S}'_\sigma(\mathbb{R}^3):=\big\{u\in[\mathscr{S}'(\mathbb{R}^3)]^3:\textrm{div}u=0\big\}$,
$M_1$, $M_3$ and $\widetilde{\mathbb{P}}$ are defined as (2.5) (2.7) and (2.8), respectively.

By the similar argument of \cite{Iwabuchi 2014},
we arrive at
\begin{eqnarray}
&&~~\big|\mathscr{F}[A_2(f^M)(t)](\xi)\big|\nonumber\\
&\geq&~~\bigg|\int^t_0\hat{T}_1(\xi,t-\tau)\sum^3_{l=1}\bigg(\delta_{1,l}-\frac{\xi_1\xi_l}{|\xi|^2}\bigg)
\sum^3_{k=1}\xi_k
[
\big(\widehat{T_{\Omega,N}(\tau)f^M}\big)_k\ast \big(\widehat{T_{\Omega,N}(\tau)f^M}\big)_l\big]d\tau\bigg|\nonumber\\
&&-\bigg|\int^t_0\hat{T}_2(\xi,t-\tau)M_2
\mathscr{F}\big[\widetilde{\mathbb{P}}\widetilde{\nabla}\cdot
\big(T_{\Omega,N}(\tau)f^M\otimes T_{\Omega,N}(\tau)f^M\big)\big](\xi)d\tau\bigg|\nonumber\\
&&-\bigg|\int^t_0\big[\hat{T}_3(\xi,t-\tau)-\hat{T}_1(\xi,t-\tau)\big]M_3
\mathscr{F}\big[\widetilde{\mathbb{P}}\widetilde{\nabla}\cdot
\big(T_{\Omega,N}(\tau)f^M\otimes T_{\Omega,N}(\tau)f^M\big)\big](\xi)d\tau\bigg|\nonumber\\
&=:&~~K_1(\xi,t)-K_2(\xi,t)-K_3(\xi,t).\setcounter{equation}{9}
\end{eqnarray}
Here we have used the fact that
$(M_1+M_3)\mathscr{F}[\widetilde{\mathbb{P}}f]
=\mathscr{F}[\widetilde{\mathbb{P}}f]$ for $f\in[\mathscr{S}'(\mathbb{R}^3)]^4$. In the following, we divide our calculations into two parts:

 \emph{\textbf{(i) Estimates for $K_1(\xi, t)$ with $\xi\in E$}}

By the definition of $T_{\Omega,N}(\cdot)$, we see
\begin{eqnarray}
K_1(\xi, t)
&\geq&\bigg|\int^t_0\hat{T}_1(\xi,t-\tau)\big(1-\frac{\xi^2_1}{|\xi|^2}\big)
\xi_1
[
\big(\widehat{T_{\Omega,N}(\tau)f^M}\big)_1\ast \big(\widehat{T_{\Omega,N}(\tau)f^M}\big)_1\big]d\tau\bigg|\nonumber\\
&&-\sum_{k=2,3}\bigg|\int^t_0\hat{T}_1(\xi,t-\tau)\big(1-\frac{\xi^2_1}{|\xi|^2}\big)
\xi_k
[
\big(\widehat{T_{\Omega,N}(\tau)f^M}\big)_k\ast \big(\widehat{T_{\Omega,N}(\tau)f^M}\big)_1\big]d\tau\bigg|\nonumber\\
&&-\sum_{k=1,2,3,l=2,3}\bigg|\int^t_0\hat{T}_1(\xi,t-\tau)\frac{\xi_1\xi_l}{|\xi|^2}
\xi_k
[
\big(\widehat{T_{\Omega,N}(\tau)f^M}\big)_k\ast \big(\widehat{T_{\Omega,N}(\tau)f^M}\big)_l\big]d\tau\bigg|\nonumber\\
&=:&K_{111}(\xi,t)-\sum_{k=2,3}K_{1k1}(\xi,t)-\sum_{k=1,2,3,l=2,3}K_{1kl}(\xi,t).
\setcounter{equation}{10}
\end{eqnarray}

\textit{On the estimate of $K_{111}(\xi,t)$ for $\xi\in E$.} By the definition of $T_{\Omega,N}$, one has
\begin{eqnarray}
K_{111}(\xi,t)
&\geq&\bigg|\int^t_0\hat{T}_1(\xi,t-\tau)\big(1-\frac{\xi^2_1}{|\xi|^2}\big)
\xi_1
[
\big(\hat{T}_{3}(\xi,\tau)M_3\widehat{f^M}\big)_1\ast \big(\hat{T}_{3}(\xi,\tau)M_3\widehat{f^M}\big)_1\big]d\tau\bigg|\nonumber\\
&&-\sum_{k=1,2}\bigg|\int^t_0\hat{T}_1(\xi,t-\tau)\big(1-\frac{\xi^2_1}{|\xi|^2}\big)
\xi_1
[
\big(\hat{T}_{k}(\xi,\tau)M_k\widehat{f^M}\big)_1\ast \big(\hat{T}_{3}(\xi,\tau)M_3\widehat{f^M}\big)_1\big]d\tau\bigg|\nonumber\\
&&-\sum_{k=1,2,3,l=1,2}\bigg|\int^t_0\hat{T}_1(\xi,t-\tau)\big(1-\frac{\xi^2_1}{|\xi|^2}\big)
\xi_1
[
\big(\hat{T}_{k}(\xi,\tau)M_k\widehat{f^M}\big)_1\ast \big(\hat{T}_{l}(\xi,\tau)M_l\widehat{f^M}\big)_1\big]d\tau\bigg|\nonumber\\
&=:&J_{133}-\sum_{k=1,2}J_{1k3}-\sum_{k=1,2,3,l=1,2}J_{1kl}.
\setcounter{equation}{11}
\end{eqnarray}

We now estimate each term. It follows from the definitions of $f^M$ that
\begin{eqnarray}
J_{133}~&=&~\bigg|2\int^t_0\hat{T}_1(\xi,t-\tau)\big(1-\frac{\xi^2_1}{|\xi|^2}\big)
\xi_1\int_{\mathbb{R}^3}
e^{-\nu|\xi-\eta|^2\tau}e^{-\nu|\eta|^2\tau}
\frac{N^2|\xi_h-\eta_h|^2}{|\xi-\eta|'^2}\frac{\xi_2-\eta_2}{|\xi-\eta|}\nonumber\\
&&~~~\times
\frac{N^2|\eta_h|^2}{|\eta|'^2}\frac{\eta_2}{|\eta|}
\frac{1}{M}\sum^{2M}_{j=M}2^{2j}\chi^+_j(\xi-\eta)\chi^-_j(\eta)d\eta d\tau\bigg|,
\nonumber
\end{eqnarray}
where we have used the support properties of $\chi^\pm_j$ and $\eta_h:=(\eta_1,\eta_2)$ for every $\eta=(\eta_1,\eta_2,\eta_3)\in\mathbb{R}^3$.

Since
$$\frac{|\xi|'}{|\xi|}(t-\tau)
\leq Nt\leq1\ \ \ \text{and}\ \ \ e^{-\frac{\nu}{N}}\leq e^{-\nu|\xi|^2(t-\tau)} $$
for $t\in(0,\frac{1}{N}]$ and $\xi\in E$,
there exists a constant $0<c<1$ such that
$$\hat{T}_1(\xi,t-\tau)\geq c,$$
and since
$$-\frac{N^4}{\Omega^4}\leq\frac{N^2|\xi_h-\eta_h|^2}{|\xi-\eta|'^2}\frac{\xi_2-\eta_2}{|\xi-\eta|}
\frac{N^2|\eta_h|^2}{|\eta|'^2}
\frac{\eta_2}{|\eta|}\leq-\frac{1}{256}$$
for $\eta\in\textrm{supp}\chi_j^-$ with $\xi-\eta\in\textrm{supp}\chi_j^+$, or
$\eta\in\textrm{supp}\chi_j^+$ with $\xi-\eta\in\textrm{supp}\chi_j^-$,
we then obtain
\begin{eqnarray}
J_{133}~&\geq&~\frac{c}{M}\sum^{2M}_{j=M}2^{2j}\int^t_0
\int_{\mathbb{R}^3}
e^{-\nu|\xi-\eta|^2\tau}e^{-\nu|\eta|^2\tau}
\chi^+_j(\xi-\eta)\chi^-_j(\eta)d\eta d\tau
\nonumber\\
&\geq&~\frac{c}{M}\sum^{2M}_{j=M}2^{2j}2^{-2j}(1-e^{-\nu4t2^{2j}})\geq c\nonumber
\end{eqnarray}
for $\xi\in E$ and $\nu^{-1}2^{-2M}\leq t\leq \frac{1}{N}$ with $M\geq \frac{1}{2}\log_2\frac{N}{\nu}$,
which yields that
$$\|J_{133}\|_{L^1(E)}\geq c.\eqno(3.12)$$

For $J_{113}$, we have for $\xi\in E$ that
\begin{eqnarray}
J_{113}~&\leq&~C\bigg|\int^t_0
\int_{\mathbb{R}^3}
e^{-\nu(|\xi-\eta|^2+|\eta|^2)\tau}\frac{\Omega^2|\xi_3-\eta_3|^2}{|\xi-\eta|'^2}\frac{N^2|\eta_h|^2}{|\eta|'^2}
\frac{\xi_2-\eta_2}{|\xi-\eta|}\frac{\eta_2}{|\eta|}\nonumber\\
&&\ \ \ \ \ \ \ \ \ \ \
\times\frac{1}{M}\sum^{2M}_{j=M}2^{2j}\chi^+_j(\xi-\eta)\chi^-_j(\eta)d\eta d\tau\bigg|\nonumber\\
&\leq&\frac{C}{M}\frac{N^2}{\Omega^2}\sum^{2M}_{j=M}2^{2j}2^{-2j}
\int_{\mathbb{R}^3}
\frac{1-e^{-\nu(|\xi-\eta|^2+|\eta|^2)t}}{\nu(|\xi-\eta|^2+|\eta|^2)}
\chi^+_j(\xi-\eta)\chi^-_j(\eta)d\eta d\tau
\nonumber\\
&\leq&\frac{C}{M}\frac{N^2}{\Omega^2}\sum^{2M}_{j=M}2^{2j}2^{-2j}2^{-2j}
\leq\frac{C N^2}{\Omega^2 M2^{2M}}.\setcounter{equation}{13}
\end{eqnarray}
Similarly,
\begin{eqnarray}
J_{111}~&\leq&~C\bigg|\int^t_0\int_{\mathbb{R}^3}
e^{-\nu(|\xi-\eta|^2+|\eta|^2)\tau}\times\frac{\Omega^2|\xi_3-\eta_3|^2}{|\xi-\eta|'^2}\frac{\xi_2-\eta_2}{|\xi-\eta|}
\frac{\Omega^2\eta_3^2}{|\eta|'^2}
\frac{\eta_2}{|\eta|}\nonumber\\
&&\ \ \ \ \ \ \ \ \ \ \times\frac{1}{M}\sum^{2M}_{j=M}2^{2j}\chi^+_j(\xi-\eta)\chi^-_j(\eta)d\eta d\tau\bigg|\nonumber\\
&\leq&~\frac{C}{M2^{4M}}.\setcounter{equation}{14}
\end{eqnarray}
Thus, noting that $J_{131}=J_{113}$, it follows from (3.13) and (3.14) that
$$\big\|J_{111}\big\|_{L^1(E)}+\|J_{113}\big\|_{L^1(E)}+\|J_{131}\big\|_{L^1(E)}\leq \frac{C N^2}{\Omega^2 M2^{2M}}. \eqno(3.15)$$

For $J_{112}$, we note that
$$|\hat{T}_2(\zeta,\tau)|\leq\frac{|\zeta|'}{|\zeta|}\tau\leq N t\eqno(3.16)$$
for all $\tau\in(0,t)$ and all $\zeta\in\mathbb{R}^3$. Hence, for $\xi\in E$ one has
\begin{eqnarray}
J_{112}&\leq&Ct\int^t_0
\int_{\mathbb{R}^3}
e^{-\nu(|\xi-\eta|^2+|\eta|^2)\tau}
\big|\widehat{f^M}(\xi-\eta)\big|~\big|\widehat{f^M}(\eta)\big|d\eta d\tau.\setcounter{equation}{17}
\end{eqnarray}
Since it holds that $\sum_{j=-j_0}^{j_0}\hat{\psi}_j(\xi)=1$ for all $\xi\in E$,
Young's inequality together with H\"{o}lder inequality, Lemma 2.1 and Lemma 2.4 ensures that
\begin{eqnarray}
\big\|J_{112}\big\|_{L^1(E)}
&\leq&Ct
\bigg\|\int^t_0
\bigg[\big[e^{-\nu|\cdot|^2\tau}
\big|\widehat{f^M}(\cdot)\big|\big]\ast\big[e^{-\nu|\cdot|^2\tau}
\big|\widehat{f^M}(\cdot)\big|\big]\bigg]d\tau\bigg\|_{L^1}\nonumber\\
&\leq&Ct
\bigg\{\sum^{j_0}_{j=-j_0}\bigg(\int^t_0\bigg\|\hat{\psi}_j
\mathscr{F}\bigg[\mathscr{F}^{-1}\big[e^{-\nu|\cdot|^2\tau}
\big|\widehat{f^M}\big|\big]^2\bigg]\bigg\|_{L^1}d\tau\bigg)^2\bigg\}^{\frac{1}{2}}\nonumber\\
&\leq&Ct
\bigg\|\mathscr{F}^{-1}\big[e^{-\nu|\cdot|^2\tau}
\big|\widehat{f^M}\big|\big]\bigg\|_{L^\frac{1}{1-\alpha}(0,\infty;\dot{FB}^{-\alpha}_{1,2})}
\bigg\|\mathscr{F}^{-1}\big[e^{-\nu|\cdot|^2\tau}
\big|\widehat{f^M}\big|\big]\bigg\|_{L^\frac{1}{1+\alpha}(0,\infty;\dot{FB}^{\alpha}_{1,2})}\nonumber\\
&\leq&Ct
\|f^M\|^2_{\dot{FB}^{-1}_{1,2}}
\leq Ct.\setcounter{equation}{18}
\end{eqnarray}
Similarly, we get
$$\|J_{121}\|_{L^1(E)}+\|J_{123}\|_{L^1(E)}+\|J_{132}\|_{L^1(E)}\leq Ct \eqno(3.19)$$
and
\begin{eqnarray}
\big\|J_{122}\big\|_{L^1(E)}
&\leq&Ct^2
\bigg\|\int^t_0
\bigg[\big[e^{-\nu|\cdot|^2\tau}
\big|\widehat{f^M}(\cdot)\big|\big]\ast\big[e^{-\nu|\cdot|^2\tau}
\big|\widehat{f^M}(\cdot)\big|\big]\bigg]d\tau\bigg\|_{L^1}\nonumber\\
&\leq&Ct^2.\setcounter{equation}{20}
\end{eqnarray}
Therefore, by (3.12), (3.15) and (3.18)-(3.20), we have
$$\|K_{111}\|_{L^1(E)}\geq c-\frac{C N^2}{\Omega^2 M2^{2M}}-C(t+t^2).\eqno(3.21)$$

\textit{On the estimate of $K_{121}(\xi,t)$ and $K_{112}(\xi,t)$ for $\xi\in E$.}
By the definition of $T_{\Omega,N}(\cdot)$, we see
\begin{eqnarray}
K_{121}(\xi,t)
&\leq&\sum_{k,l=1,2,3}\bigg|\int^t_0
[
\big(\hat{T}_{k}(\xi,\tau)M_k\widehat{f^M}\big)_2\ast \big(\hat{T}_{l}(\xi,\tau)M_l\widehat{f^M}\big)_1\big]d\tau\bigg|\nonumber\\
&=:&\sum_{k,l=1,2,3}L_{1kl}.
\setcounter{equation}{22}
\end{eqnarray}
The similar process for getting (3.13) gives
\begin{eqnarray}
L_{111}&\leq&\bigg|-2\int^t_0\int_{\mathbb{R}^3}
e^{-\nu(|\xi-\eta|^2+|\eta|^2)\tau}\frac{\Omega^2(\xi_3-\eta_3)^2}{|\xi-\eta|'^2}
\frac{\xi_1-\eta_1}{|\xi-\eta|}\frac{\Omega^2\eta_3^2}{|\eta|'^2}
\frac{\eta_2}{|\eta|}\nonumber\\
&&\ \ \ \ \ \ \ \ \ \
\times\frac{1}{M}\sum^{2M}_{j=M}2^{2j}\chi^+_j(\xi-\eta)\chi^-_j(\eta)d\eta d\tau\bigg|\nonumber\\
&\leq&\frac{C}{M2^{5M}}.\setcounter{equation}{23}
\end{eqnarray}
\begin{eqnarray}
L_{113}=L_{131}&\leq&\bigg|2\int^t_0\int_{\mathbb{R}^3}
e^{-\nu(|\xi-\eta|^2+|\eta|^2)\tau}
\frac{N^2|\xi_h-\eta_h|^2}{|\xi-\eta|'^2}\frac{\xi_1-\eta_1}{|\xi-\eta|}
\frac{\Omega^2\eta_3^2}{|\eta|'^2}
\frac{\eta_2}{|\eta|}\nonumber\\
&&~~~~~~~~
\times\frac{1}{M}\sum^{2M}_{j=M}2^{2j}\chi^+_j(\xi-\eta)\chi^-_j(\eta)d\eta d\tau\bigg|\nonumber\\
&\leq&\frac{CN^2}{\Omega^2M2^{3M}},\setcounter{equation}{24}
\end{eqnarray}
and
\begin{eqnarray}
L_{133}&\leq&\bigg|2\int^t_0\int_{\mathbb{R}^3}
e^{-\nu(|\xi-\eta|^2+|\eta|^2)\tau}
\frac{N^2|\xi_h-\eta_h|^2}{|\xi-\eta|'^2}\frac{\xi_1-\eta_1}{|\xi-\eta|}
\frac{N^2|\eta_h|^2}{|\eta|'^2}
\frac{\eta_2}{|\eta|}\nonumber\\
&&~~~~~~~~
\times\frac{1}{M}\sum^{2M}_{j=M}2^{2j}\chi^+_j(\xi-\eta)\chi^-_j(\eta)d\eta d\tau\bigg|\nonumber\\
&\leq&\frac{CN^4}{\Omega^4M2^{M}},\setcounter{equation}{25}
\end{eqnarray}
which yield that
$$\big\|\sum_{k,l=1,3}L_{1kl}\big\|_{L^1(E)}\leq\frac{CN^4}{\Omega^4M2^{M}}.\eqno(3.26)$$
And the similar arguments for getting (3.18) and (3.20) give rise to
\begin{eqnarray}
\bigg\|\sum_{l=1,2,3}L_{12l}+\sum_{k=1,3}L_{1k2}\bigg\|_{L^1(E)}
&\leq&C(t+t^2)
\bigg\|\int^t_0
\bigg[\big[e^{-\nu|\cdot|^2\tau}
\big|\widehat{f^M}(\cdot)\big|\big]\ast\big[e^{-\nu|\cdot|^2\tau}
\big|\widehat{f^M}(\cdot)\big|\big]\bigg]d\tau\bigg\|_{L^1}\nonumber\\
&\leq&C(t+t^2).\setcounter{equation}{27}
\end{eqnarray}

Therefore, it follows from (3.26) and (3.27) that
$$\|K_{121}\|_{L^1(E)}\leq \frac{CN^4}{\Omega^4M2^{M}}+C(t+t^2).\eqno(3.28)$$
Similarly,
$$\|K_{112}\|_{L^1(E)}\leq \frac{CN^4}{\Omega^4M2^{M}}+C(t+t^2).\eqno(3.29)$$

\vskip 3mm

\textit{On the estimate of $K_{122}(\xi,t)$ for $\xi\in E$}.
By the definition of $T_{\Omega,N}(\cdot)$, we see
\begin{eqnarray}
K_{122}(\xi,t)
&\leq&\sum_{k,l=1,2,3}\bigg|\int^t_0
[
\big(\hat{T}_{k}(\xi,\tau)M_k\widehat{f^M}\big)_2\ast \big(\hat{T}_{l}(\xi,\tau)M_l\widehat{f^M}\big)_2\big]d\tau\bigg|\nonumber\\
&=:&\sum_{k,l=1,2,3}J_{2kl}.
\setcounter{equation}{30}
\end{eqnarray}
Similar calculations for getting (3.13) lead to
\begin{eqnarray}
\sum_{k,l=1,3}J_{2kl} &\leq& \bigg|-2\int^t_0\int_{\mathbb{R}^3}
e^{-\nu(|\xi-\eta|^2+|\eta|^2)\tau}\bigg\{
\frac{N^2|\xi_h-\eta_h|^2}{|\xi-\eta|'^2}\frac{N^2|\eta_h|^2}{|\eta|'^2}\nonumber\\
&&~~~~~~~~
+2\frac{N^2|\xi_h-\eta_h|^2}{|\xi-\eta|'^2}\frac{\Omega^2\eta_3^2}{|\eta|'^2}
+\frac{\Omega^2(\xi_3-\eta_3)^2}{|\xi-\eta|'^2}\frac{\Omega^2\eta_3^2}{|\eta|'^2}
\bigg\}
\frac{\xi_1-\eta_1}{|\xi-\eta|}\frac{\eta_1}{|\eta|}\nonumber\\
&&~~~~~~~~
\times\frac{1}{M}\sum^{2M}_{j=M}2^{2j}\chi^+_j(\xi-\eta)\chi^-_j(\eta)d\eta d\tau\bigg|\nonumber\\
&\leq&\frac{CN^4}{\Omega^4M2^{2M}},
\end{eqnarray}
which yields that
$$\big\|\sum_{k,l=1,3}J_{2kl}\big\|_{L^1(E)}\leq\frac{CN^4}{\Omega^4M2^{2M}}.\eqno(3.32)$$
And similar to (3.27), one has
\begin{eqnarray}
\bigg\|\sum_{l=1,2,3}J_{22l}+\sum_{k=1,3}J_{2k2}\bigg\|_{L^1(E)}
&\leq&C(t+t^2).\setcounter{equation}{33}
\end{eqnarray}

Therefore, it follows from (3.32) and (3.33) that
$$\|K_{122}\|_{L^1(E)}\leq \frac{CN^4}{\Omega^4M2^{2M}}+C(t+t^2).\eqno(3.34)$$

\textit{On the estimates of $K_{131}(\xi,t)$, $K_{113}(\xi,t)$, $K_{123}(\xi,t)$, $K_{132}(\xi,t)$ and $K_{133}(\xi,t)$ for $\xi\in E$.}
By the definition of $T_{\Omega,N}(\cdot)$, together with the fact that $(M_1\widehat{f^M})_3=0$ and $(M_3\hat{v})_3=0$ for all $v\in[\mathscr{S}'(\mathbb{R}^3)]^4$, one has
\begin{eqnarray}
K_{131}(\xi,t)
&\leq&\sum_{l=1,2,3}\bigg|\int^t_0
[
\big(\hat{T}_{2}(\xi,\tau)M_2\widehat{f^M}\big)_3\ast \big(\hat{T}_{l}(\xi,\tau)M_l\widehat{f^M}\big)_1\big]d\tau\bigg|,\setcounter{equation}{35}
\end{eqnarray}
\begin{eqnarray}
K_{123}(\xi,t)
&\leq&\sum_{k=1,2,3}\bigg|\int^t_0
[
\big(\hat{T}_{k}(\xi,\tau)M_k\widehat{f^M}\big)_2\ast \big(\hat{T}_{2}(\xi,\tau)M_2\widehat{f^M}\big)_3\big]d\tau\bigg|\setcounter{equation}{36}
\end{eqnarray}
and
\begin{eqnarray}
K_{133}(\xi,t)
&\leq&\bigg|\int^t_0
[
\big(\hat{T}_{2}(\xi,\tau)M_2\widehat{f^M}\big)_3\ast \big(\hat{T}_{2}(\xi,\tau)M_2\widehat{f^M}\big)_3\big]d\tau\bigg|.\setcounter{equation}{37}
\end{eqnarray}
Similar to (3.27), we obtain
$$
\|K_{131}\|_{L^1(E)}+\|K_{123}\|_{L^1(E)}+\|K_{133}\|_{L^1(E)}
\leq C(t+t^2).\eqno(3.38)
$$
And similarly,
$$
\|K_{113}\|_{L^1(E)}+\|K_{132}\|_{L^1(E)}
\leq C(t+t^2).\eqno(3.39)
$$

Summing up (3.21), (3.28), (3.29), (3.34), (3.38) and (3.39), we get
$$\|K_1\|_{L^1{(E)}}\geq c-\frac{CN^4}{\Omega^4M2^{M}}-C(t+t^2)\eqno(3.40)$$
for $-\nu2^{-2M}\leq t\leq \frac{1}{N}$ with $M\geq\frac{1}{2}\log_2\frac{N}{\nu}$.

\vskip 3mm

\emph{\textbf{(ii) Estimates for $K_2(\xi, t)$ and $K_3(\xi, t)$ with $\xi\in E$}}

Noticing that
$$|\hat{T}_2(\zeta,\tau)|\leq\frac{|\zeta|'}{|\zeta|}\tau\leq N t$$
and
$$|\hat{T}_3(\zeta,\tau)-\hat{T}_1(\zeta,\tau)|
\leq1-\cos\bigg(\frac{|\zeta|'}{|\zeta|}\tau\bigg)
\leq\frac{1}{2}\frac{|\zeta|'^2}{|\zeta|^2}\tau^2\leq\frac{1}{2}N^2 t^2$$
for all $\tau\in(0,t)$ and all $\zeta\in\mathbb{R}^3$, it is easy to observe for $\xi\in E$
that
$$
K_{2}(\xi, t)+K_{3}(\xi, t)\leq C(t+t^2)\int^t_0
\int_{\mathbb{R}^3}
e^{-\nu(|\xi-\eta|^2+|\eta|^2)\tau}
\big|\widehat{f^M}(\xi-\eta)\big|~\big|\widehat{f^M}(\eta)\big|d\eta d\tau.
$$
The same process for getting (3.18) results in
\begin{eqnarray}
\big\|K_{2}+K_{3}\big\|_{L^1(E)}
&\leq&C(t+t^2).\setcounter{equation}{41}
\end{eqnarray}

Finally, summing up (3.40) and (3.41), we see that there exist $0<c$ and $C>1$ such that
$$\|A_2(f^M)(t)\|_{\dot{FB}^{-1}_{1,r}}\geq c-\frac{CN^4}{\Omega^4M2^{M}}-C(t+t^2)$$
for $\nu^{-1}2^{-2M}\leq t\leq \frac{1}{N}$ with $M\geq\frac{1}{2}\log_2\frac{N}{\nu}$.
In fact, there holds
$$\|A_2(f^M)(t)\|_{\dot{FB}^{-1}_{1,r}}\geq \frac{c}{3}$$
for all $M\in\mathbb{N}$ and $t>0$ with $M\geq\max\{\frac{3CN^4}{\Omega^4c},  \frac{1}{2}\log_2\frac{6C}{c\nu}, \frac{1}{2}\log_2\frac{N}{\nu},\frac{1}{2}\log_2\frac{1}{\nu}\}$ and $\nu^{-1}2^{-2M}\leq t\leq\min\{\frac{c}{6C}, \frac{1}{N}, 1\}$, which implies (3.3).
This completes the proof of Theorem 1.6.
\rightline{$\Box$}
\vskip 2mm

\setlength{\baselineskip}{6pt}

\bibliographystyle{abbrv}

\end{document}